\documentclass[twoside]{article}

\oddsidemargin=\evensidemargin
\catcode`\@=11
\long\def\@makefntext#1{
\protect\noindent \hbox to 3.2pt {\hskip-.9pt  
$^{{\eightrm\@thefnmark}}$\hfil}#1\hfill}		%CAN BE USED 
	
\def\ps@myheadings{\let\@mkboth\@gobbletwo		%SIZE OF R/H NOS.
\def\@oddhead{\hbox{}
\rightmark\hfil\eightrm\thepage}   
\def\@oddfoot{}\def\@evenhead{\eightrm\thepage\hfil
\leftmark\hbox{}}\def\@evenfoot{}
\def\sectionmark##1{}\def\subsectionmark##1{}}

\catcode`@=11		      		     %SIZE OF OPENING PAGE NO.
\def\ps@plain{\let\@mkboth\@gobbletwo
     \def\@oddhead{}\def\@oddfoot{\eightrm\hfil\thepage
     \hfil}\def\@evenhead{}\let\@evenfoot\@oddfoot}

\newcounter{sectionc}\newcounter{subsectionc}\newcounter{subsubsectionc}
\renewcommand{\section}[1] {\vspace{12pt}\addtocounter{sectionc}{1} 
\setcounter{subsectionc}{0}\setcounter{subsubsectionc}{0}\noindent 
	{\tenbf\thesectionc. #1}\par\vspace{5pt}}
\renewcommand{\subsection}[1] 
{\vspace{12pt}\addtocounter{subsectionc}{1} 
	\setcounter{subsubsectionc}{0}\noindent 
	{\bf\thesectionc.\thesubsectionc. 
	{\kern1pt \bfit #1}}\par\vspace{5pt}}
\renewcommand{\subsubsection}[1] {\vspace{12pt}
	\addtocounter{subsubsectionc}{1}
	\noindent
	{\tenrm\thesectionc.\thesubsectionc.\thesubsubsectionc.	{\kern1pt 
	\it #1}}\par\vspace{5pt}}
\newcommand{\nonumsection}[1] {\vspace{12pt}\noindent{\tenbf #1}
	\par\vspace{5pt}}

\newcommand{\textlineskip}{\baselineskip=13pt}
\newcommand{\smalllineskip}{\baselineskip=10pt}

\newcommand{\copyrightheading}[1]
	{\vspace*{-2.5cm}\smalllineskip{\flushleft
	{\footnotesize Journal of Knot Theory and Its Ramifications #1}\\
   	{\footnotesize \copyright\kern2pt World Scientific 
         Publishing Company}\\
         }}

\def\abstracts#1{{
	\footnotesize\baselineskip=10pt
	\noindent{\bf Abstract.} #1
	\par}} 

\def\keywords#1{{ 
	\footnotesize\baselineskip=10pt
	\noindent{\bf Keywords.} #1
	\par}}
	
\def\AMS#1{{ 
	\footnotesize\baselineskip=10pt
	\noindent{\bf Mathematics Subject Classification (2000).} #1
	\par}}

\renewenvironment{thebibliography}[1]
	{\frenchspacing

	 \ninerm\baselineskip=11pt
	 \begin{list}{[\arabic{enumi}]}
	{\usecounter{enumi}\setlength{\parsep}{0pt}
	 \setlength{\leftmargin 13.7pt}{\rightmargin 0pt} %[1--9] ITEMS
	 \setlength{\itemsep}{0pt} \settowidth
	{\labelwidth}{[#1]}\sloppy}}{\end{list}}

\newcounter{itemlistc}
\newcounter{Mlistc}
\newcounter{romanlistc}

\newenvironment{itemlist}
    	{\setcounter{itemlistc}{0}
	 \begin{list}{-}
	{\usecounter{itemlistc}
	 \setlength{\parsep}{0pt}
	\setlength{\topsep}{4pt}
 	\setlength{\itemsep}{0pt}}}{\end{list}}
	 
\newenvironment{Mlist}
    	{\setcounter{Mlistc}{0}
	 \begin{list}{(M\arabic{Mlistc})}
	{\usecounter{Mlistc}
	 \setlength{\parsep}{0pt}
	 \setlength{\topsep}{4pt}
	 \setlength{\itemsep}{0pt}}}{\end{list}}

\newenvironment{romanlist}
	{\setcounter{romanlistc}{0}
	 \begin{list}{$($\roman{romanlistc}$)$}
	{\usecounter{romanlistc}
	 \setlength{\topsep}{4pt}
	 \setlength{\parsep}{0pt}
	 \setlength{\itemsep}{0pt}}}{\end{list}}

\newcommand{\fcaption}[1]{
        \refstepcounter{figure}
        \setbox\@tempboxa = \hbox{\footnotesize Fig.~\thefigure. #1}
        \ifdim \wd\@tempboxa > 5in
           {\begin{center}
        \parbox{5in}{\footnotesize\smalllineskip Fig.~\thefigure. #1}
            \end{center}}
        \else
             {\begin{center}
             {\footnotesize Fig.~\thefigure. #1}
              \end{center}}
        \fi}

\def\pmb#1{\setbox0=\hbox{#1}
	\kern-.025em\copy0\kern-\wd0
	\kern.05em\copy0\kern-\wd0
	\kern-.025em\raise.0433em\box0}

\def\fnt#1#2{\footnotetext{\kern-.3em
	{$^{\mbox{\scriptsize #1}}$}{#2}}}

\def\fpage#1{\begingroup
\voffset=.3in
\thispagestyle{empty}\begin{table}[b]\centerline{\footnotesize #1}
	\end{table}\endgroup}

\def\runninghead#1#2{\pagestyle{myheadings}
\markboth{{\protect\footnotesize{\centerline{#1}}}}
{\hfill{\protect\footnotesize{\centerline{#2}}}}}

\font\tenrm=cmr10
 
\font\tenbf=cmbx10
\font\bfit=cmbxti10 at 10pt
\font\ninerm=cmr9
\font\nineit=cmti9
\font\ninebf=cmbx9
\font\eightrm=cmr8

\def\qed{\hbox{${\vcenter{\vbox{			%HOLLOW SQUARE
   \hrule height 0.4pt\hbox{\vrule width 0.4pt height 6pt
   \kern5pt\vrule width 0.4pt}\hrule height 0.4pt}}}$}}

\usepackage{epsfig}
\usepackage{floatfig}
\usepackage{amssymb}
\newcommand{\Z}{{\mathbb Z}}
\newcommand{\Ho}{\hbox{H}}
\newcommand{\sgn}{\hbox{sgn}}
\newcommand{\rk}{\hbox{rk}\,}
\newcommand{\N}{{\cal N}}
\newcommand{\cqfd}{\hfill \qed \medskip}
\begin{document}

\initfloatingfigs

\setlength{\textheight}{7.7truein}  %for 2nd page onwards  

\runninghead{D. Cimasoni}{Conway potential function}

\normalsize\textlineskip
\thispagestyle{empty}
\setcounter{page}{1}

\fpage{1}
\centerline{\bf A geometric construction of the Conway potential function}
\baselineskip=13pt
\vspace*{0.37truein}
\vspace{0.5cm}															  
\noindent David Cimasoni
\baselineskip=12pt

\vspace*{0.5truein}

\abstracts{We give a geometric construction of the multivariable Conway potential function for
colored links. In the case of a single color, it is Kauffman's definition of the Conway polynomial in terms of a Seifert matrix.}
\vspace*{10pt}
\AMS{57M25.}
\vspace*{10pt}
\keywords{Conway potential function, Conway polynomial, colored link, C-complex.}

\vspace{1cm}															

\vspace*{1pt}\textlineskip
\section{Introduction}
\vspace*{-0.5pt}

A {\it colored link\/} is an oriented link $L=L_1\cup\dots\cup L_\mu$ in $S^3$ together with  a surjective map
$\sigma$ assigning to each component $L_i$ a color $\sigma(i)\in\{1,\dots,n\}$. Two colored links $(L,\sigma)$ and $(L',\sigma')$
are {\it isotopic\/} if there exists an ambient isotopy from $L$ to $L'$ preserving the color and orientation of each component.

Given a colored link $(L,\sigma)$ with exterior $X$, we have $\Ho_1X=\bigoplus_{i=1}^\mu\Z t_i$, where $t_1,\dots,t_\mu$ denote the oriented
meridians of $L$. The Hurewicz morphism $\pi_1X\longrightarrow\Ho_1X$ composed with
$\bigoplus_{i=1}^\mu\Z t_i\longrightarrow\bigoplus_{i=1}^n\Z t_i\,,\;\;t_i\mapsto t_{\sigma(i)}$, determines a regular
$\Z^n$-covering $\widehat{X}^\sigma\longrightarrow X$. The homology of $\widehat{X}^\sigma$ is endowed with a natural
structure of module over $\Lambda_n=\Z[t_1,t_1^{-1},\dots,t_n,t_n^{-1}]$. The $\Lambda_n$-module $\Ho_1\widehat{X}^\sigma$ is
called the {\it Alexander module\/} of $(L,\sigma)$: it is an invariant of the colored link.
The greatest commun divisor of the first elementary ideal of this module
is called the {\it multivariable Alexander polynomial\/} of $(L,\sigma)$.
It is denoted by $\Delta_L^\sigma(t_1,\dots,t_n)$.
Note that $\Delta_L^\sigma$ is only defined modulo the units of $\Lambda_n$, that is, up to multiplication by
$\pm t_1^{\nu_1}\cdots t_n^{\nu_n}$.

If $n=1$, $\Delta_L^\sigma(t)$ is nothing but the Alexander polynomial of the unordered oriented link $L$, as defined by {\sc James Alexander} \cite{Alex}.
On the other hand, if $n=\mu$ and $\sigma=\hbox{id}$, $\Delta_L^\sigma(t_1,\dots,t_\mu)$ is denoted by $\Delta_L(t_1,\dots,t_\mu)$: it is the Alexander
polynomial of the ordered oriented link $L$, defined by {\sc Ralph Fox} \cite{Fox}. Note that in general,
$$
\Delta^{\sigma}_L(t_1,\dots,t_n)\,\;\dot{=}\,\;
\cases{(t_1-1)\Delta_L(t_1,\dots,t_1) & if $n=1$ and $\mu>1$;\cr
	\Delta_L(t_{\sigma(1)},\dots,t_{\sigma(\mu)}) & else,}
$$
where $\dot{=}$ stands for the equality modulo the units of $\Lambda_n$.

In 1970, {\sc John Conway} \cite{Con} introduced a new invariant of links called the {\it potential function}.
Given a colored link $(L,\sigma)$, its potential function is a well-defined rational function
$\nabla^\sigma_L(t_1,\dots,t_n)$ which satisfies
$$
\nabla^\sigma_L(t_1,\dots,t_n)\,\;\dot{=}\,\cases{\frac{1}{t_1-t_1^{-1}}\,\Delta^\sigma_L(t_1^2) & if $n=1$; \cr           
            	\Delta^\sigma_L(t_1^2,\dots,t_n^2) & if $n>1$.\cr}
$$
Thus, this invariant is just the Alexander polynomial without the ambiguity
concerning multiplication by units of $\Lambda_n$. This might seem a minor improvement. However,
the potential function has a very remarkable new property: it can be computed directly from a link diagram using skein relations.
For example, suppose that $(L_+,\sigma)$ $(L_-,\sigma)$ and $(L_0,\sigma_0)$ are given by colored diagrams related by a single change as below, where
$i$ stands for the color of the components.

\begin{figure}[h]
\begin{center}
\epsfig{figure=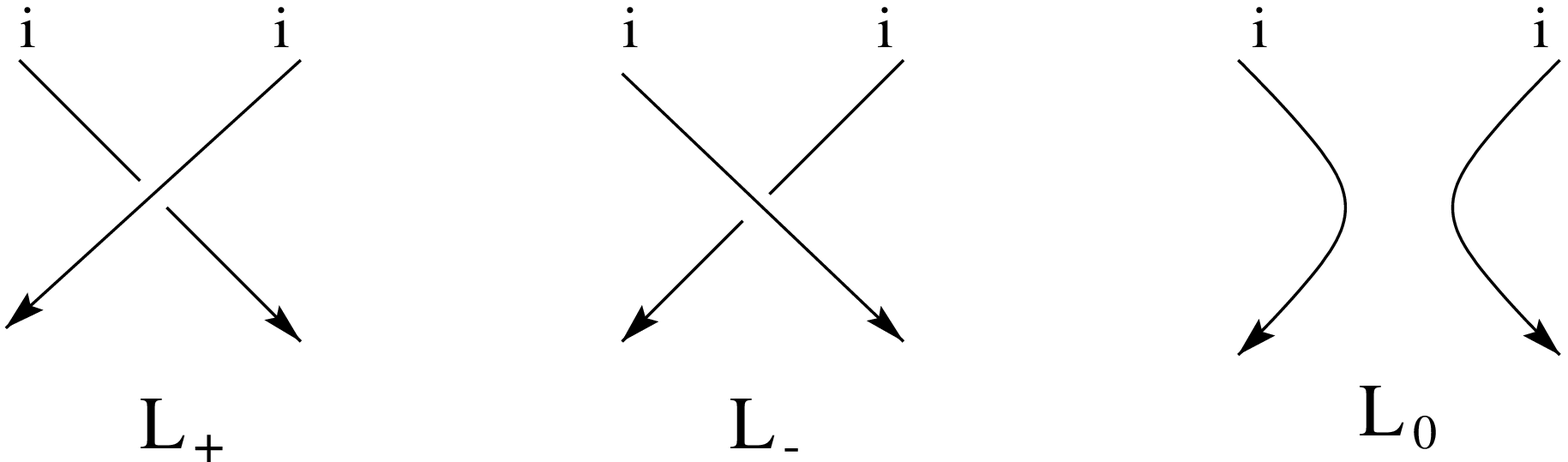,height=1.8cm}
\end{center}
\end{figure}
\noindent Then, their potential functions satisfy the equality
$$
\nabla^\sigma_{L_+}-\nabla^\sigma_{L_-}=(t_i-t_i^{-1})\cdot \nabla^{\sigma_0}_{L_0}.
$$
Similarly, if $L_{++}$, $L_{--}$ and $L_{00}$ differ by the following local operation,

\begin{figure}[h]
\begin{center}
\epsfig{figure=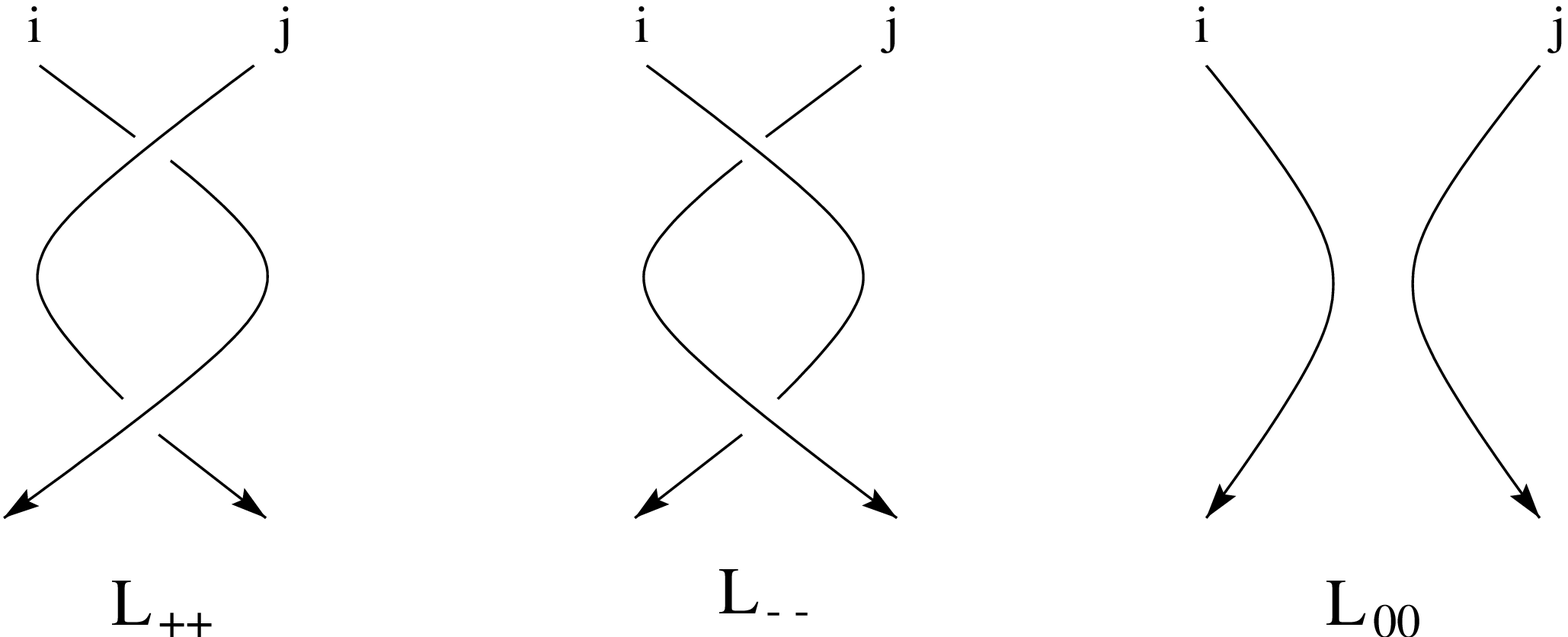,height=2.7cm}
\end{center}
\end{figure}
\noindent then we have the equality
$$
\nabla^\sigma_{L_{++}}+\nabla^\sigma_{L_{--}}=(t_it_j+t_i^{-1}t_j^{-1})\cdot \nabla^\sigma_{L_{00}}.
$$
Thus, Conway pointed out a preferred representative of the Alexander polynomial, and gave a very easy method for computing it.
Unfortunately, his paper contains neither a precise definition of the potential function, nor a proof of its unicity. Quoting Conway:
``We have not found a satisfactory {\em explanation\/} of these identities.(...) It seems plain that much work remains to be done in the field.''

As a particular case of his potential function, Conway defined what he called the {\it reduced polynomial\/} $D_L(t)$ of a non-ordered link $L$ by
$$
D_L(t)=(t-t^{-1})\cdot\nabla^\sigma_L(t),
$$
where $\sigma$ is the coloring map with a single color.
This invariant was later called the {\it Conway polynomial} of $L$; we will use this terminology.
Following from the properties of the general potential function, we have:
\begin{romanlist}
\item{$D_L$ is an invariant of the non-ordered oriented link $L$;}
\item{$D_O(t)=1$, where {\it O} stands for the trivial knot;}
\item{$D_{L_+}(t)-D_{L_-}(t)=(t-t^{-1})\cdot D_{L_0}(t)$, where $L_+,L_-$ and $L_0$ are related by a single crossing change as below.}
\end{romanlist}

\begin{figure}[h]
\begin{center}
\epsfig{figure=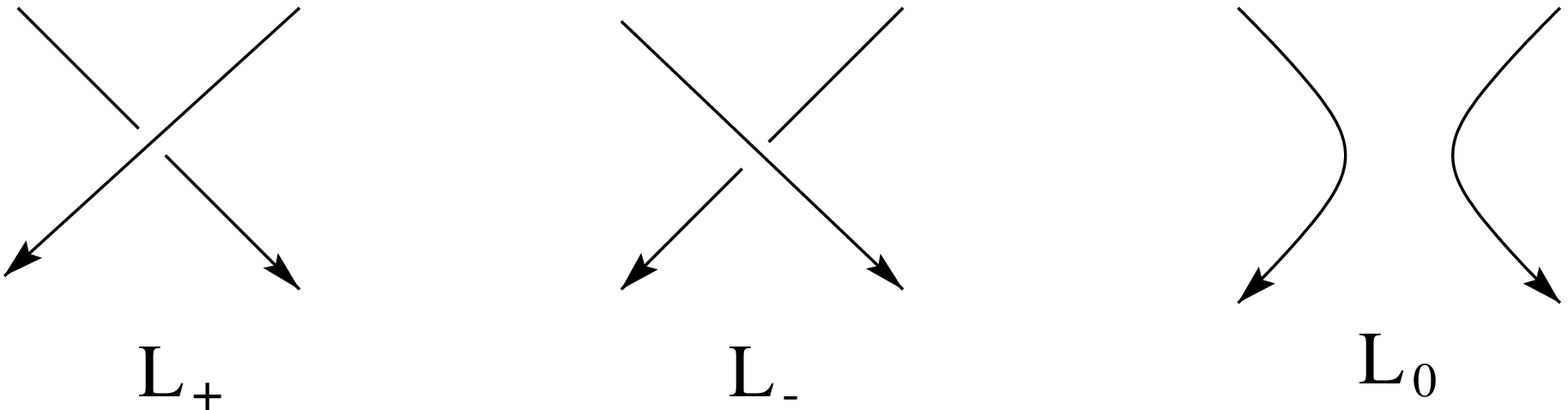,height=1.8cm}
\end{center}
\end{figure}
But again, the paper \cite{Con} contains no proof of the constitency of this system of axioms: there is no explicit model for $D_L$. On the other hand,
it is very easy to check that if such a model exists, it is unique.

In 1981, {\sc Louis Kauffman} \cite{Kauf} found a very simple geometric construction of the Conway polynomial $D_L(t)$, namely
$$
D_L(t)=\det\left(-tA+t^{-1}A^T\right),
$$
where $A$ is any Seifert matrix of the link $L$ and $A^T$ is the transpose of $A$. This result gave the requested model for $D_L$; it also
provided the `explanation' of the first skein relation. 

Finally, in 1983, {\sc Richard Hartley} \cite{Har} gave a definition of the
multivariable potential function $\nabla^\sigma_L$. Quoting his introduction:
``Kauffman showed how to define what may be called the reduced potential function of a link in terms of a Seifert matrix. This reduced potential
function is an $L$-polynomial in one variable. However, the potential function is essentially a function of several variables, and I can see no
way of generalising Kauffman's method to obtain the full potential function. Quite a different approach is therefore indicated.'' Indeed,  Hartley's
definition is obtained by normalizing the Alexander matrix given by the Wirtinger presentation via Fox free differential calculus: it is an algebraic
construction, whose relation to Kauffman's model is not obvious at all.

In this paper, we give a geometric construction of the multivariable potential function which generalizes Kauffman's definition.
It gives the relation between Hartley and Kauffman's models of the Conway polynomial, and provides an `explanation' of Conway's skein
formulas. It should be pointed out that the fundamental idea of this construction, that is, the use of `C-complexes' for the computation of
Alexander invariants of links, is due to {\sc Daryl Cooper} \cite{Cooper1,Cooper2}.

\section{The potential function of a colored link}

Let us denote by $\nabla_L$ the potential function of an ordered oriented link $L$ as defined by Hartley \cite{Har}; this corresponds to the case
of a colored link $(L,\sigma)$ with $n=\mu$ and $\sigma=\hbox{id}$. The potential function $\nabla_L^\sigma$ of an arbitrary colored link
is given by
$$
\nabla^{\sigma}_L(t_1,\dots,t_n)=\nabla_L(t_{\sigma(1)},\dots,t_{\sigma(\mu)}).
$$
The aim of this section is to present a geometric construction of $\nabla^{\sigma}_L$. This requires the use of C-complexes
as defined by Cooper \cite{Cooper1,Cooper2}.

\begin{figure}[Htb]
\begin{center}
\epsfig{figure=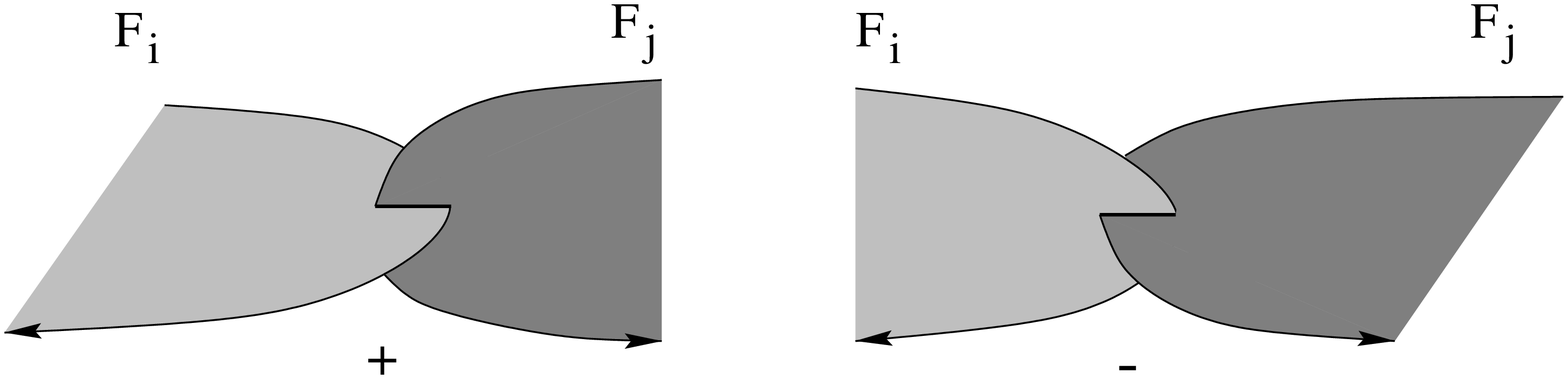,height=2.8cm}
\fcaption{A positive clasp, and a negative clasp.}
\label{fig:clasps}
\end{center}
\end{figure}

\medskip
\noindent {\bf Definition.} {\it A C-complex for a colored link $(L,\sigma)$ is a union $F=F_1\cup\dots\cup F_n$
of compact oriented PL-embedded surfaces in $S^3$ which satisfies the following properties:
\begin{romanlist}
\item{F is connected;}
\item{for all $i$, $\partial F_i$ is equal to the sublink of $L$ of color $i$;}
\item{for all $i\neq j$, $F_i\cap F_j$ is either empty, or a disjoint union of clasps (see Figure~\ref{fig:clasps});}
\item{for all $i,j,k$ pairwise distinct, $F_i\cap F_j\cap F_k$ is empty.}
\end{romanlist}}
\medskip

C-complexes are a natural generalization of Seifert surfaces; indeed, if $n=1$, a C-complex for $(L,\sigma)$ is simply
a Seifert surface for $L$. Now, we need to define the Seifert forms associated to a C-complex. 
Let us fix a C-complex $F$ and a map $\epsilon\colon\{1,\dots,n\}\to\{\pm 1\}$. A $1$-cycle $x$ in $F$ is called a {\it loop\/}
if it has the following behaviour near a clasp (whenever it crosses one).

\begin{figure}[h]
\begin{center}
\epsfig{figure=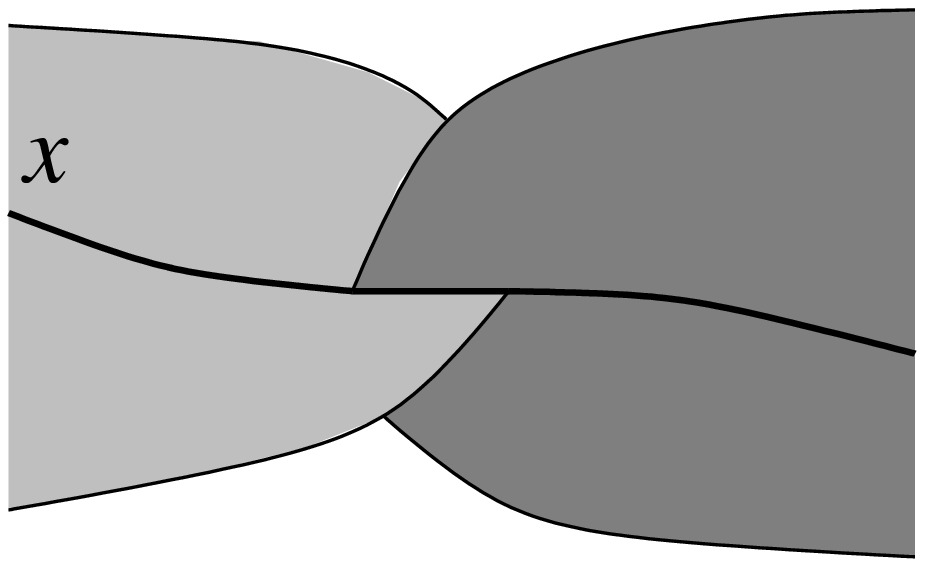,height=2.3cm}
\end{center}
\end{figure}
Given any element in $\Ho_1F$, it is possible to represent it by a loop; therefore, we can define a bilinear
form $\alpha^\epsilon\colon\Ho_1F\times\Ho_1F\to\Z$ by 
$$
\alpha^\epsilon([x],[y])=\ell k(x^\epsilon,y),
$$
where $\ell k$ denotes the linking number, $x$ is a loop, and $x^\epsilon$ is obtained by pushing $x$ off $F$ in the $\epsilon(i)-$normal direction
off $F_i$ for all $i=1,\dots,n$. In the neighborhood of a clasp in $F_i\cap F_j$, a loop $x$ can be pushed off $F$ in four different ways,
which correspond to the four possible values of $(\epsilon(i),\epsilon(j))$. This is illustrated by the following cross-section of a clasp
in $F_i\cap F_j$.
\begin{figure}[h]
\begin{center}
\epsfig{figure=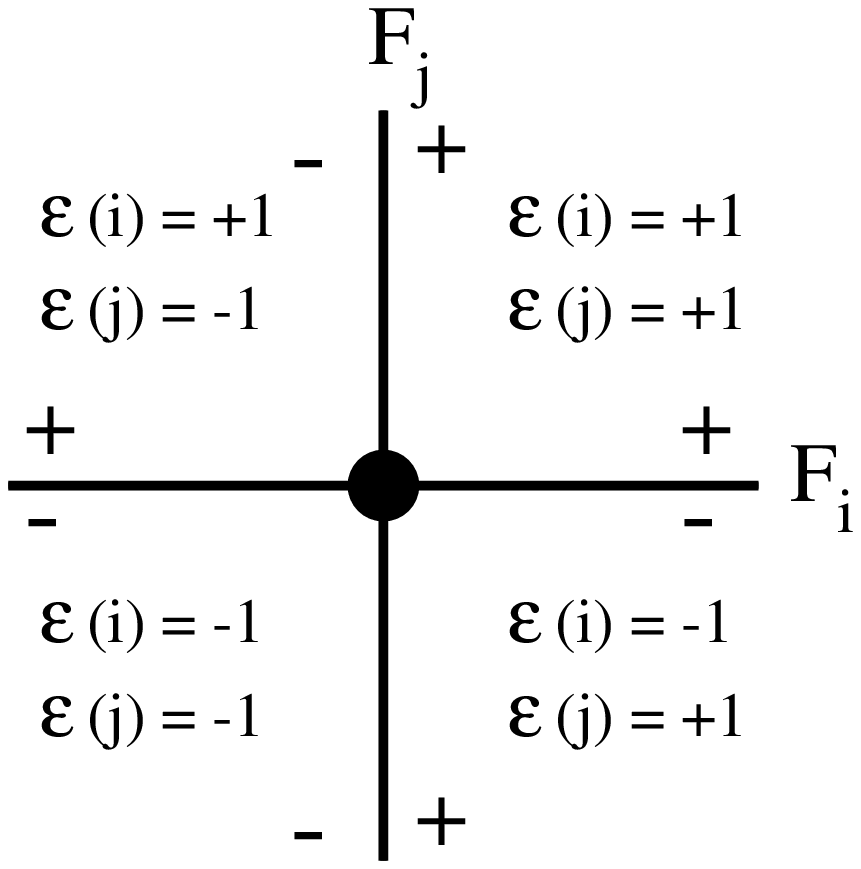,height=2.7cm}
\end{center}
\end{figure}

Of course, if $n=1$, $\alpha^-$ is the usual Seifert form.
Let us fix a basis of $\Ho_1F$ and denote by $A_F^\epsilon$ the matrix of $\alpha^\epsilon$. Note that for all
$\epsilon$, $A_F^{-\epsilon}=\left(A_F^\epsilon\right)^T$. We are now in position to state our main result.
\medskip

\noindent {\bf Theorem.} {\it The Conway potential function of a colored link $(L,\sigma)$ is given by
$$
\nabla_L^\sigma(t_1,\dots,t_n)=\sgn(F)\prod_{i=1}^n\left(t_i-t_i^{-1}\right)^{\chi(F\setminus F_i)-1}\det\left(-A_F\right),
$$
where $F$ is any C-complex for $(L,\sigma)$, and where
\begin{itemlist}
\item{sgn$(F)$ is equal to the product of the signs of all the clasps in $F$ (see Figure~\ref{fig:clasps});}
\item{$A_F=\sum_\epsilon\epsilon(1)\cdots\epsilon(n)\cdot t_1^{\epsilon(1)}\cdots t_n^{\epsilon(n)}\cdot A_F^\epsilon$,
the sum being on the $2^n$ maps $\epsilon\colon\{1,\dots,n\}\to\{\pm 1\}$.}
\end{itemlist}
}\medskip

\noindent {\bf Example~1.} If $n=1$, $F$ is a Seifert surface for the link $L$ and $A_F^-$ is nothing but a usual Seifert matrix $A$. Since
$A^+_F=A^T$, we have
$$
\nabla^\sigma_L(t)=\frac{1}{t-t^{-1}}\cdot\det\left(-tA^T+t^{-1}A\right),
$$
giving Kauffman's construction of the Conway polynomial
$$
D_L(t)=(t-t^{-1})\cdot\nabla^\sigma_L(t)=\det\left(-tA+t^{-1}A^T\right).
$$

\noindent {\bf Example~2.} If $n=2$, let us denote by $x,y$ the variables and by $F=F_x\cup F_y$ a C-complex for $(L,\sigma)$. Let us also
note $A_F^{--}=A$ and $A_F^{-+}=B$. We then have the formula 
$$
\nabla_L^\sigma(x,y)=\sgn(F)(x-x^{-1})^{\chi(F_y)-1}(y-y^{-1})^{\chi(F_x)-1}\det\left(-xyA^T+xy^{-1}B^T+x^{-1}yB-x^{-1}y^{-1}A\right).
$$
Compare Cooper \cite[Corollary 2.2]{Cooper2}.
\medskip

\begin{floatingfigure}{4cm}
\mbox{\epsfig{figure=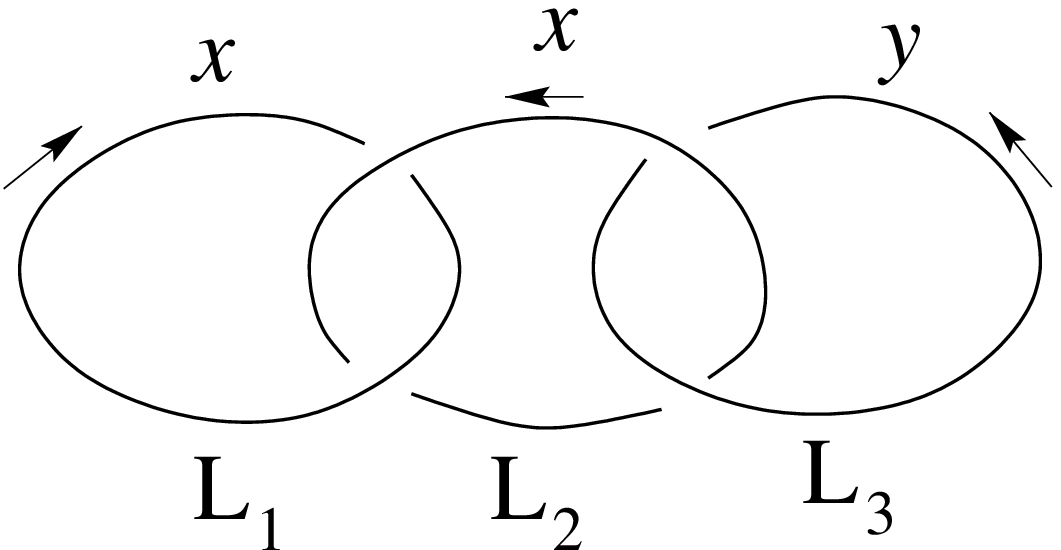,width=3.5cm}}
\end{floatingfigure}
\noindent {\bf Example~3.} Consider the link $L=L_1\cup L_2\cup L_3$ with colors $\sigma(1)=\sigma(2)=x$ and $\sigma(3)=y$ given by the following diagram.
We have two different ways to compute the potential function $\nabla_L^\sigma(x,y)$ using our construction. The first one is to consider a C-complex for the
colored link $(L,\sigma)$, and to compute $\nabla_L^\sigma(x,y)$ directly; the second one is to calculate $\nabla_L(t_1,t_2,t_3)$ using a C-complex for
the ordered link $L$, and to use the equality $\nabla_L^\sigma(x,y)=\nabla_L(x,x,y)$. Here are C-complexes $F=F_x\cup F_y$ for $(L,\sigma)$ and
$F'=F_1\cup F_2\cup F_3$ for the ordered link $L$.

\begin{figure}[h]
\begin{center}
\epsfig{figure=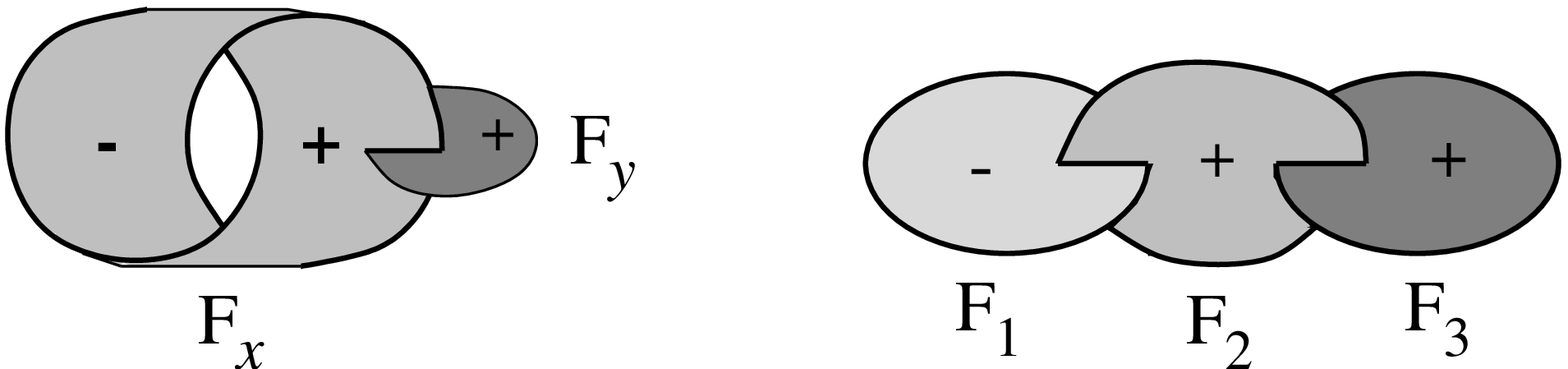,height=2cm}
\end{center}
\end{figure}
\noindent In the obvious basis of $\Ho_1F$, we have $A_F^{--}=A_F^{-+}=(-1)$; therefore, the theorem gives
$$
\nabla_L^\sigma(x,y)=(x-x^{-1})^0\cdot(y-y^{-1})^{-1}\cdot(xy-x^{-1}y-xy^{-1}+x^{-1}y^{-1})=x-x^{-1}.
$$
On the other hand, $F'$ is contractible, giving
$$
\nabla_L(t_1,t_2,t_3)=(t_1-t_1^{-1})^0\cdot(t_2-t_2^{-1})^1\cdot(t_3-t_3^{-1})^0=t_2-t_2^{-1},
$$
and we get the same result $\nabla_L^\sigma(x,y)=\nabla_L(x,x,y)=x-x^{-1}$.
 
\section{Proof of the theorem}

Let us recall very briefly Kauffman's argument \cite{Kauf}: given an oriented link $L$, he sets $\Omega_L(t)=\det\left(-tA+t^{-1}A^T\right)$, where
$A$ is any Seifert surface for $L$.
The first step is to check that $\Omega_L$ is an isotopy invariant of $L$. This is achieved using a well known result:
two Seifert surfaces for ambient isotopic links can be obtained from each other by a finite number of ambient isotopies and
handle attachments. We are left with the easy proof that $\Omega_L$ is unchanged by attaching a handle to $F$. Then, Kauffman
proves that $\Omega_L$ satisfies the equality
$$
\Omega_{L_+}(t)-\Omega_{L_-}(t)=(t-t^{-1})\cdot \Omega_{L_0}(t),
$$
where $L_+,L_-$ and $L_0$ are related by a single crossing change as described in the introduction.
This relation, along with the value of $\Omega$ for the trivial knot, determines $\Omega$ for any link. Since this relation
was among the formulas announced by Conway for his polynomial $D$, and since $\Omega$ and $D$ have the same value on
the trivial knot, they are the same invariant.

We will follow the same procedure for the general case. First, we define $\Omega^\sigma_L$ to be the right-hand side of the equality
in the theorem. To show that $\Omega_L^\sigma$ is a well-defined isotopy invariant of colored links, we will generalize Cooper's Isotopy Lemma
\cite[Lemma~3.2]{Cooper2}. A theorem of Murakami \cite{Mur} states that the potential function $\nabla_L^\sigma$ is determined by several
skein relations. Therefore, we will just need to check that $\Omega^\sigma_L$ satisfies these properties in order to have $\Omega=\nabla$.\footnote{We
should also mention that Turaev \cite{Tur} gave another characterization of Conway's multivariable potential function. One could also
use his set of axioms to prove that $\Omega_L$ is indeed the potential function.}
\bigskip

\begin{figure}[Htb]
\begin{center}
\epsfig{figure=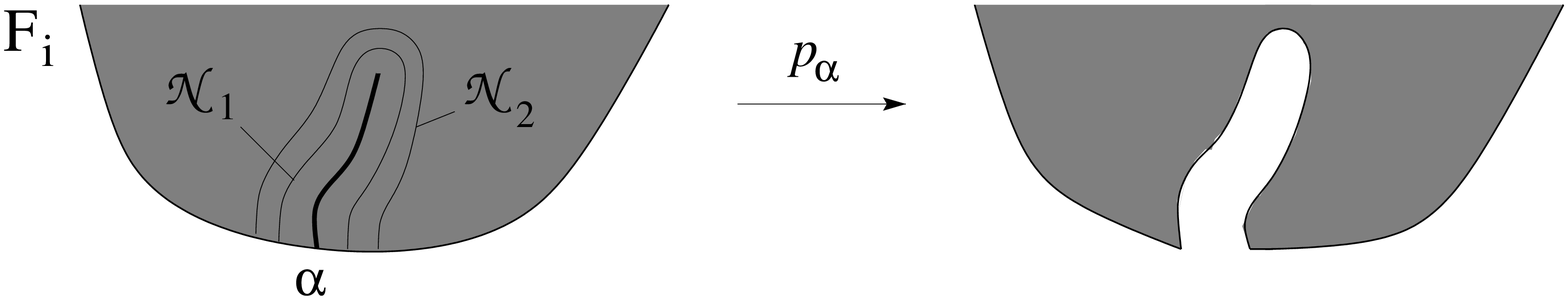,height=1.8cm}
\fcaption{Push in $F_i$ along an arc $\alpha$.}
\label{fig:push}
\end{center}
\end{figure}
As a preliminary, let us define three transformations of C-complexes; the first two were introduced by Cooper.
Given a C-complex $F$ and an arc $\alpha\colon[0,1]\to F_i$ with $\alpha^{-1}(\partial F)=\{0\}$, a {\it push in $F_i$ along\/} $\alpha$ is
an embedding $p_\alpha\colon F\to F$ defined as follows: choose $\N_1$, $\N_2$ two closed regular neighborhoods of $\alpha$ in $F_i$ meeting
$\partial F$ regularly such that $\N_1\subset\hbox{int}\,\N_2$ and $\N_2\cap(\partial F\setminus\partial F_i)$ is empty. Then, $p_\alpha$ restricted to
$F\setminus\hbox{int}\,\N_2$ is the identity, and $p_\alpha$ maps $\N_2$ homeomorphically onto $\N_2\setminus\hbox{int}\,\N_1$,
as in Figure~\ref{fig:push}. Note that it is not allowed to push a boundary component through another one.
\begin{figure}[Htb]
\begin{center}
\epsfig{figure=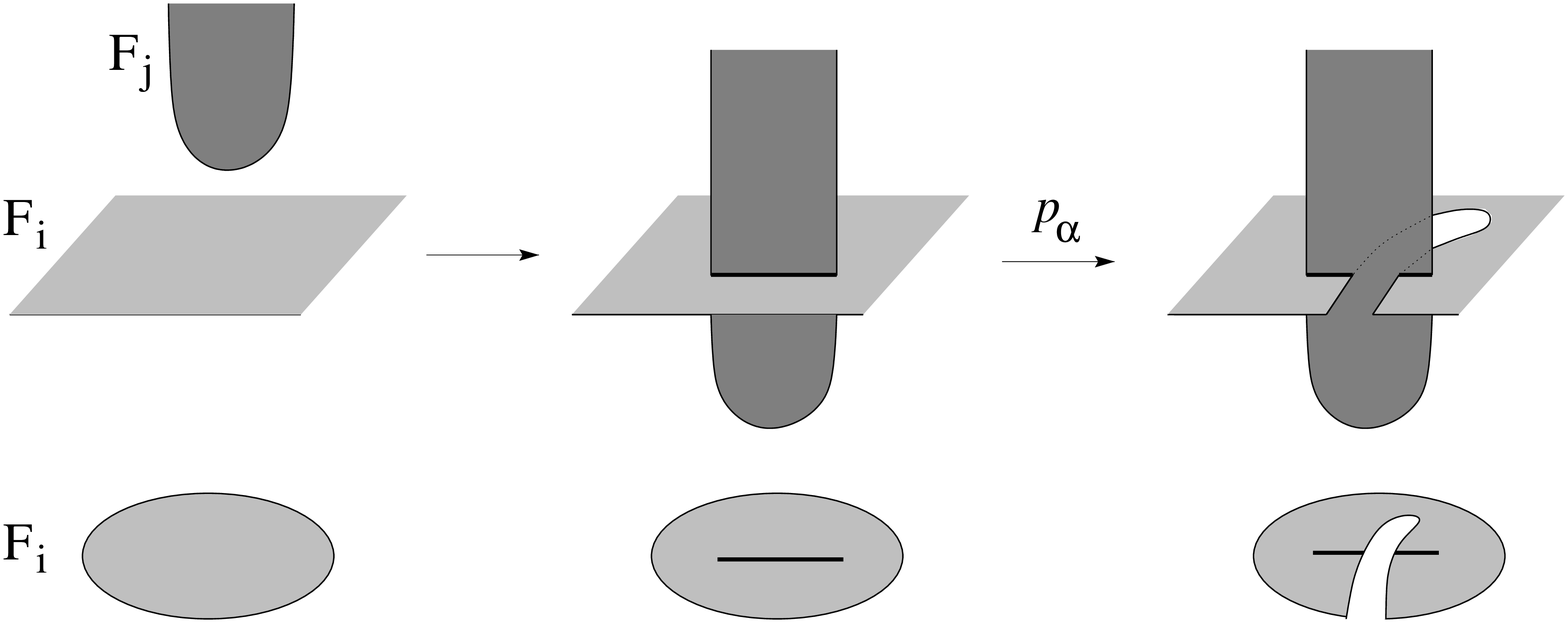,height=3.8cm}
\fcaption{Add a ribbon intersection, and push along an arc to convert it into two clasps.}
\label{fig:M1}
\end{center}
\end{figure}

One can also transform a C-complex by {\it adding a ribbon intersection\/}, as described in the left part of Figure~\ref{fig:M1}. The result is no longer
a C-complex, but it can be transformed into one by pushing in $F_i$ along an arc through the ribbon intersection: this converts the ribbon intersection into a pair
of clasps (see the right part of Figure~\ref{fig:M1}). The lower part of this figure shows the same transformations as the upper part, with the
Seifert surface $F_i$ represented as a disc.

\begin{figure}[Htb]
\begin{center}
\epsfig{figure=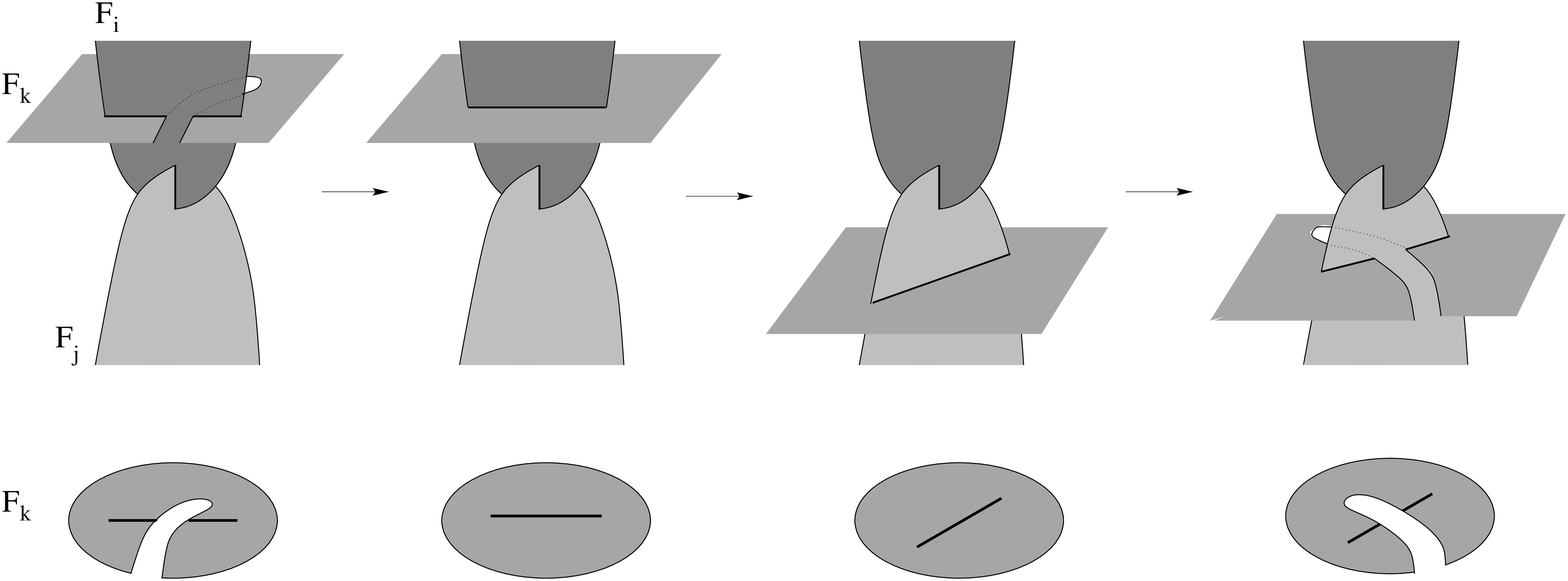,height=4.5cm}
\fcaption{Pass through a clasp.}
\label{fig:pass}
\end{center}
\end{figure}

Finally, the transformation illustrated in Figure~\ref{fig:pass} is called to {\it pass through a clasp\/}. Again, the lower part of this figure depicts
the same transformation, with a disc representing the surface $F_k$.
Note that these three transformations move the link $L$ only up to an ambient isotopy. Therefore, they can be understood as keeping $L=\partial F$ fixed.

We are now ready to start the proof of the theorem. Let $(L,\sigma)$ be a colored link in $S^3$.
\medskip

\noindent{\bf Lemma~1.} {\it Let $F_1,\dots,F_n$ be a collection of Seifert surfaces for the colored link $(L,\sigma)$ 
(i.e. for all $i$, $\partial F_i$ is equal to the sublink of $L$ of color $i$). Then, each $F_i$ can be isotoped keeping its boundary fixed
to give a C-complex for $(L,\sigma)$.}\medskip

\begin{figure}[Htb]
\begin{center}
\epsfig{figure=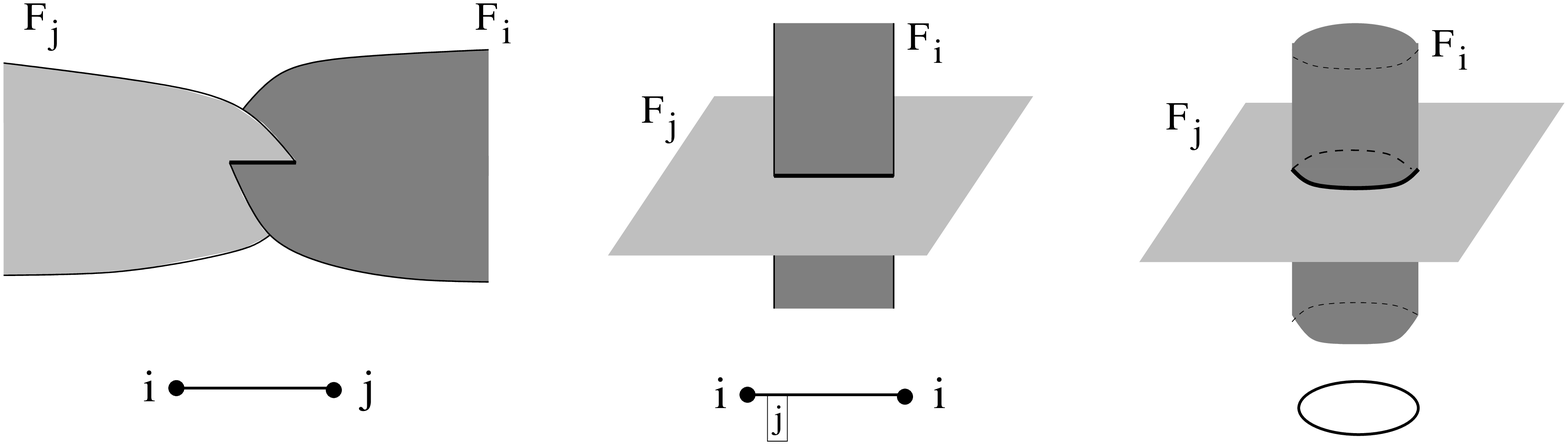,height=3.8cm}
\fcaption{The three types of intersections $F_i\cap F_j$: a clasp, a ribbon and a circle.}
\label{fig:int}
\end{center}
\end{figure}

\noindent{\it Proof.} Although it is an easy generalization of Cooper's \cite[Lemma~3.1]{Cooper2}, it is worth giving the proof in some details.
Let $F_1,\dots,F_n$ be a collection of Seifert surfaces for $(L,\sigma)$. By isotopies, it may be assumed that they
intersect transversally; therefore:
\begin{itemlist}
\item{for all $i\neq j$, $F_i\cap F_j$ is a finite union of intervals ({\it clasps\/} or {\it ribbons\/}) and {\it circles} (see Figure~\ref{fig:int});}
\item{for all $i,j,k$ pairwise distinct, $F_i\cap F_j\cap F_k$ is a finite number of points (called {\it triple points\/});}
\item{every quadruple intersection is empty.}
\end{itemlist}

By pushing along arcs, it is possible to remove every circle intersection as follows. First, push along an arc on $F_1$ through one of the outermost
circles on $F_1$, avoiding the triple points: it will turn this circle intersection into a ribbon intersection. This method allows to remove
every circle intersection on $F_1$. Use the same procedure for $F_2,\dots,F_n$, turning every circle into a ribbon.

Now that all the intersections $F_i\cap F_j$ are finite unions of intervals, it is easy to remove all the triple points: pick an interval, number its
possible triple points starting at one end and finishing at the other, and remove them in this order using the isotopy illustrated below.

\begin{figure}[h]
\begin{center}
\epsfig{figure=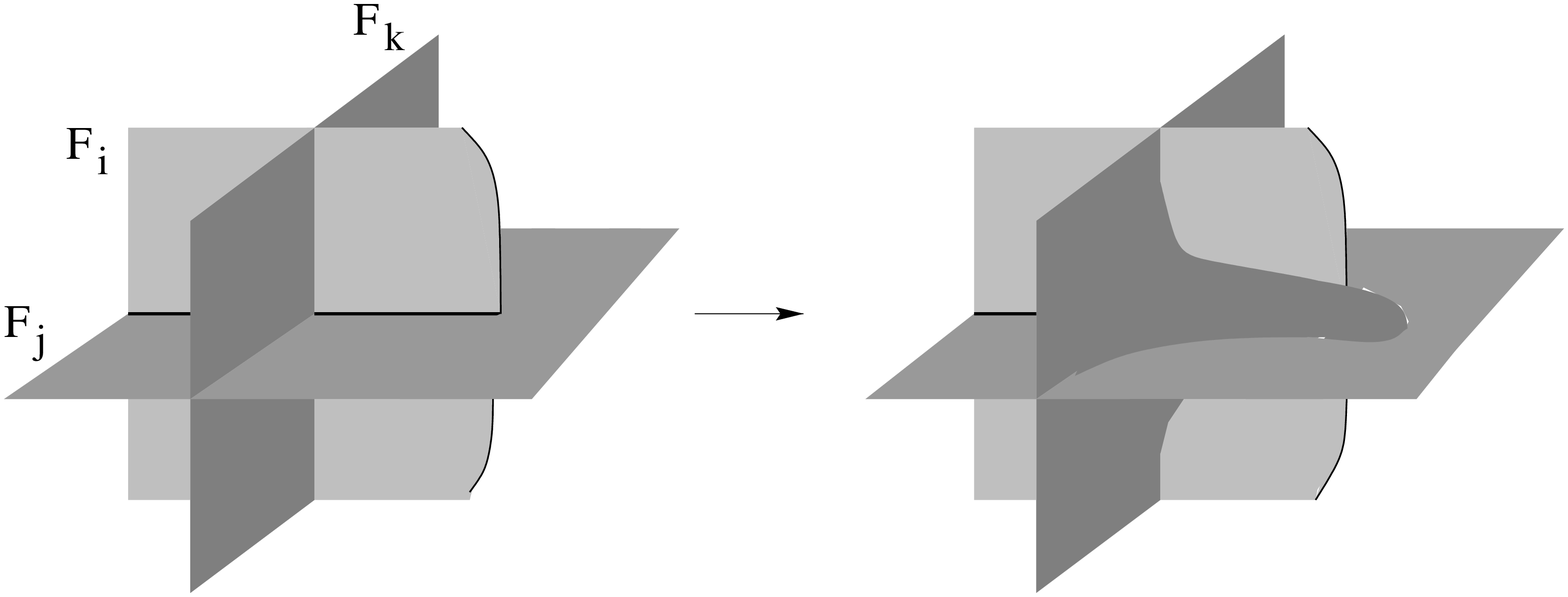,height=3.2cm}
\end{center}
\end{figure}
By pushing along arcs, it is possible to transform all the ribbon intersections into pairs of clasps, as on the right part of Figure~\ref{fig:M1}.
We now have a union $F=F_1\cup\dots\cup F_n$ which satisfies conditions (ii), (iii) and (iv) in the 
definition of a C-complex. By adding ribbon intersections and pushing along arcs (as in Figure~\ref{fig:M1}), we can make $F$ connected. \cqfd

Given $F$ a C-complex for $(L,\sigma)$, let us define
$$
\Omega_F(t_1,\dots,t_n)=\sgn(F)\prod_{i=1}^n\left(t_i-t_i^{-1}\right)^{\chi(F\setminus F_i)-1}\det\left(-A_F\right),
$$
where sgn$(F)$ is equal to the product of the signs of all the clasps in $F$, and where
$A_F=\sum_\epsilon\epsilon(1)\cdots\epsilon(n)\cdot t_1^{\epsilon(1)}\cdots t_n^{\epsilon(n)}\cdot A_F^\epsilon$,
the sum being on the $2^n$ different maps $\epsilon\colon\{1,\dots,n\}\to\{\pm 1\}$. To prove that $\Omega_F$ does not depend on the choice of
$F$, we need two lemmas. The first one is well known (see e.g. \cite[Theorem~8.2]{Lick}).
\medskip

\noindent{\bf Lemma~2.} {\it Two Seifert surfaces for a fixed link $L$ can be transformed into each other by a finite number of the following
operations and their inverses:
\begin{itemlist}
\item{ambient isotopy keeping $L$ fixed;}
\item{addition of a handle. \cqfd}
\end{itemlist}}

The second is a generalization of Cooper's Isotopy Lemma \cite[Lemma~3.2]{Cooper2}, which corresponds to the case $\mu=n=2$.
\medskip

\noindent{\bf Lemma~3.} {\it Let $F=F_1\cup\dots\cup F_n$ and $F'=F'_1\cup\dots\cup F'_n$ be two C-complexes for a fixed colored link $(L,\sigma)$.
Suppose that, for all $i$, $F_i$ is ambient isotopic to $F_i'$ keeping their boundary fixed. Then, $F$ and $F'$ can be transformed into each other by a
finite number of the following operations and their inverses:
\begin{Mlist}
\setcounter{Mlistc}{-1}
\item{ambient isotopy keeping $L$ fixed;}
\item{addition of a ribbon intersection, followed by a push along an arc through this intersection;}
\item{pass through a clasp.}
\end{Mlist}}

\noindent{\it Proof.} By \cite[Lemma~5.2]{Cooper1}, there are other ambient isotopies $F_i\sim F_i'\;(rel\; \partial F_i)$ such that the surfaces $F_1,\dots,F_n$
remain transverse throughout the isotopies except at a finite number of critical points, each of which occurs at a different time.
Let us carry out these new isotopies one after the other: the result is an ambient isotopy $F\sim F'\;(rel\;L)$ with finitely many critical points.
Each critical point corresponds to a change in the {\it singularity\/} of $F$, that is, in the intersections $F_i\cap F_j$ and
$F_i\cap F_j\cap F_k$ (the isotopies may be chosen avoiding quadruple intersections). We need to list all the possible transformations of these intersections,
and convert them into combinations of the operations stated in the lemma.

Recall that the connected components of $F_i\cap F_j$ can be of the following three types: a clasp, a ribbon, or a circle (see Figure~\ref{fig:int}).
The first step is to show that the isotopy may be chosen such that no circle intersection appears. The idea is to transform the isotopy by pushing along
wandering arcs: a push in $F_i$ along a wandering arc is an embedding $p_\alpha\colon F_i\times I\to F_i\times I$ such that $p_\alpha$ restricted to
$F_i\times\{t\}$ is a push along an arc $\alpha_t\colon  I\to F_i$.
A circle intersection in $F_i\cap F_j$ can be avoided as follows: push in $F_i$ along a wandering arc to the circle just before its appearance,
keep the arc breaking the circle during the whole `lifetime' of the circle, and remove the arc after the `death' of the circle.
Therefore, it may be assumed that all the intersections $F_i\cap F_j$ are clasps and ribbons.

Given any $i\neq j$, the number of clasps in $F_i\cap F_j$ is equal modulo $2$ to the linking number $\ell k(\partial F_i,\partial F_j)$.
Therefore, the parity of this number must be preserved.
Considering critical points that are not triple points, we are left with the following possible transformations, where $*$ indicates the critical point.

\begin{figure}[h]
\begin{center}
\epsfig{figure=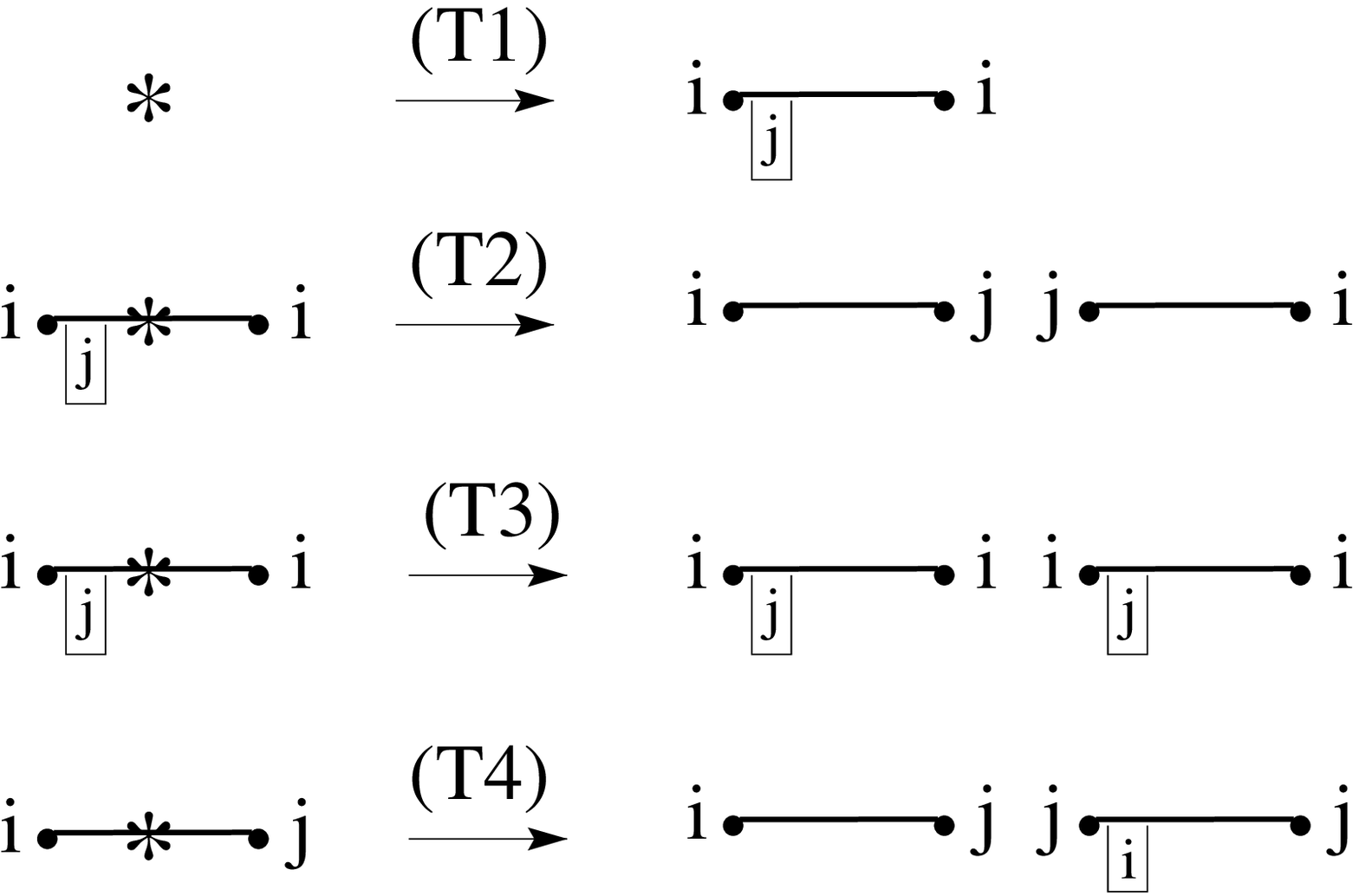,height=3.5cm}
\end{center}
\end{figure}
Transformation $T1$ is the addition of a ribbon intersection, while $T2$ is a push along an arc through a ribbon intersection. Clearly, $T3$ can be
converted into $T2 + T4 - T2$, where $-T2$ denotes the inverse of $T2$. Finally, $T4$ may be turned into $T1$ and $T2$ via a push in
$F_i$ along a wandering arc.

Let us now consider critical points that are triple points, and denote by $T5$ the move illustrated by the second arrow in Figure~\ref{fig:pass}. Here is the
`movie' of $T5$, where the asterisk denotes the critical points of `birth' and `death' of the  triple point.

\begin{figure}[h]
\begin{center}
\epsfig{figure=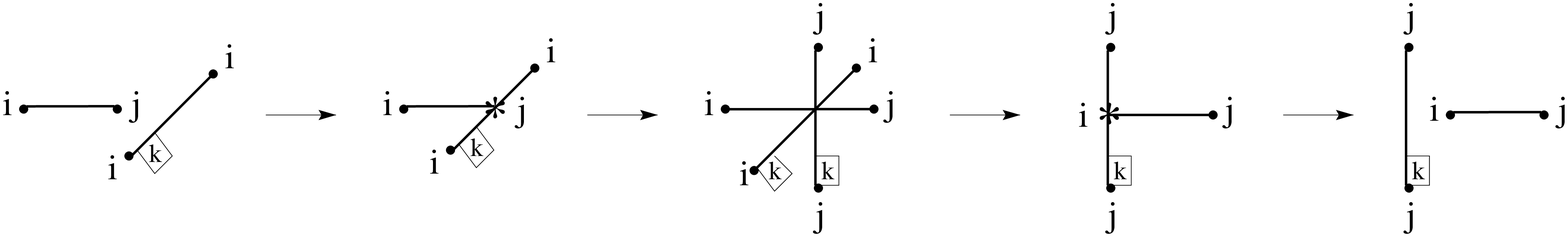,height=2cm}
\end{center}
\end{figure}
Of course, there are other transformations involving triple points; let us show that they can be converted into combinations of $T1$, $T2$
and $T5$. By the beginning of the proof, we just need to convert such a move into combinations of $T5$ and of transformations involving no triple points.
Let us consider a transformation of a C-complex $F$ involving a triple point $x\in F_1\cap F_2\cap F_3$, and let us focus on the connected component
$S$ of the singularity of $F$ containing $x$. By pushing along wandering arcs, it may be assumed that the only critical points on $S$ during the whole
`lifetime' of $x$ are the two critical points of `birth' and `death' of $x$. The idea is now to list all the possible configurations of $S$ (as in the
middle of the `movie' of $T5$), and for every such $S$, to give all the different moves involving this $S$. Finally, we get seven possible configurations
of $S$, and a dozen of (essentially) different moves. Then, it is an easy but tedious exercice to check that all these moves can be converted into combinations
of $T5$ and of transformations that do not involve triple points.

Hence, two C-complexes $F$ and $F'$ as in the statement of the lemma can be transformed into each other by a finite number of the operations
$M0$, $T1$, $T2$, $T5$ and their inverses. Since $F$ and $F'$ are C-complexes, the transformations must start with $T1$ or $-T2$ and finish with $T2$ or $-T1$.
Now, every combination of $\pm T1, \pm T2$ and $\pm T5$ starting with $T1$ or $-T2$ and finishing with $T2$ or $-T1$ can be written as a combination of
$\pm M1=\pm(T1+T2)$ and $\pm M2=\pm(-T2+T5+T2)$. This settles the proof. \cqfd
 
\noindent{\bf Lemma~4.} {\it Let $F$ and $F'$ be two C-complexes for isotopic colored links $(L,\sigma)$ and $(L',\sigma')$.
Then, $\Omega_F$ and $\Omega_{F'}$ are equal.}\medskip

\noindent{\it Proof.} Via an ambient isotopy (which does not change $\Omega_F$), we can assume that $\partial F=L=L'=\partial F'$. By Lemma~2,
$F_i$ and $F_i'$ are related by ambient isotopies (keeping $\partial F_i=\partial F'_i$ fixed) and addition of handles. These handle additions are
performed along arcs $\alpha$ embedded in $S^3\setminus F_i$; such an arc is isotopic (in $S^3\setminus F_i$) to an arc $\alpha'$ embedded in
$S^3\setminus F$. In other words, a handle attachment on $F_i$ can be performed avoiding $F\setminus F_i$. Let us call $M3$ the addition of a handle on $F_i$
avoiding the rest of $F$. Now, for every ambient isotopy between $F_i$ and $F_i'$, we can apply Lemma~3. Therefore, we just need to check that
$\Omega_F$ remains unchanged if $F$ undergoes one of the moves $M1$, $M2$ or $M3$.

\noindent{\bf (M1)} Let $F'$ be a C-complex obtained from $F$ via the move $M1$.
Clearly, $\Ho_1F'\simeq\Z\beta\oplus\Z\gamma\oplus\Ho_1F$, with $\beta$ and $\gamma$ as follows. 

\begin{figure}[h]
\begin{center}
\epsfig{figure=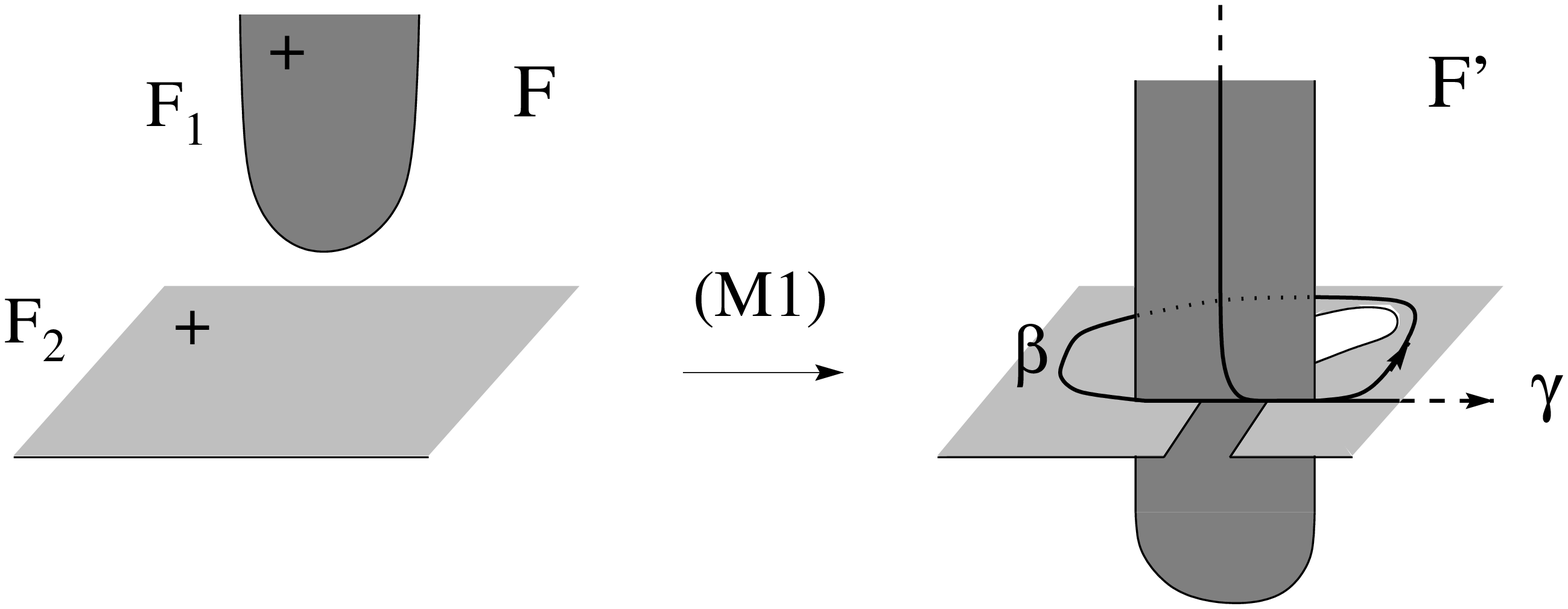,height=3.3cm}
\end{center}
\end{figure}
Among the new clasps created, one is positive and the other one is negative; therefore, $\sgn(F')=-\sgn(F)$. On the other hand,
$$
\chi(F'\setminus F'_i)=\cases{\chi(F\setminus F_i)& if $i=1,2$;\cr \chi(F\setminus F_i)-2& if $i>2$.}
$$
Since
$$
A_{F'}^\epsilon=\bordermatrix{&\beta&\gamma& \cr
			\beta&0&\delta(\epsilon)&0\cr
			\gamma&\delta(-\epsilon)&\ast&\ast\cr
			&0&\ast&A_F^\epsilon}\;,
\quad\hbox{with}\quad\delta(\epsilon)=\cases{-1& if $\epsilon(1)=\epsilon(2)=+1$;\cr\phantom{-}0& else,}
$$
it follows
$$
-A_{F'}=\pmatrix{0&t_1t_2\prod_{i>2}(t_i-t_i^{-1})&0\cr
		t_1^{-1}t_2^{-1}\prod_{i>2}(t_i-t_i^{-1})&\ast&\ast\cr
		0&\ast&-A_F}.
$$
Therefore,
\begin{eqnarray*}
\Omega_{F'}(t_1,\dots,t_n)&=&\sgn(F')\prod_i(t_i-t_i^{-1})^{\chi(F'\setminus F'_i)-1}\det(-A_{F'})\\
		&=&-\sgn(F)\frac{\prod_i(t_i-t_i^{-1})^{\chi(F\setminus F_i)-1}}{\prod_{i>2}(t_i-t_i^{-1})^2}(-1)\prod_{i>2}(t_i-t_i^{-1})^2\det(-A_{F})\\
		&=&\Omega_F(t_1,\dots,t_n).
\end{eqnarray*}

\noindent{\bf (M2)} Let $F$ and $F'$ be C-complexes as illustrated below.

\begin{figure}[h]
\begin{center}
\epsfig{figure=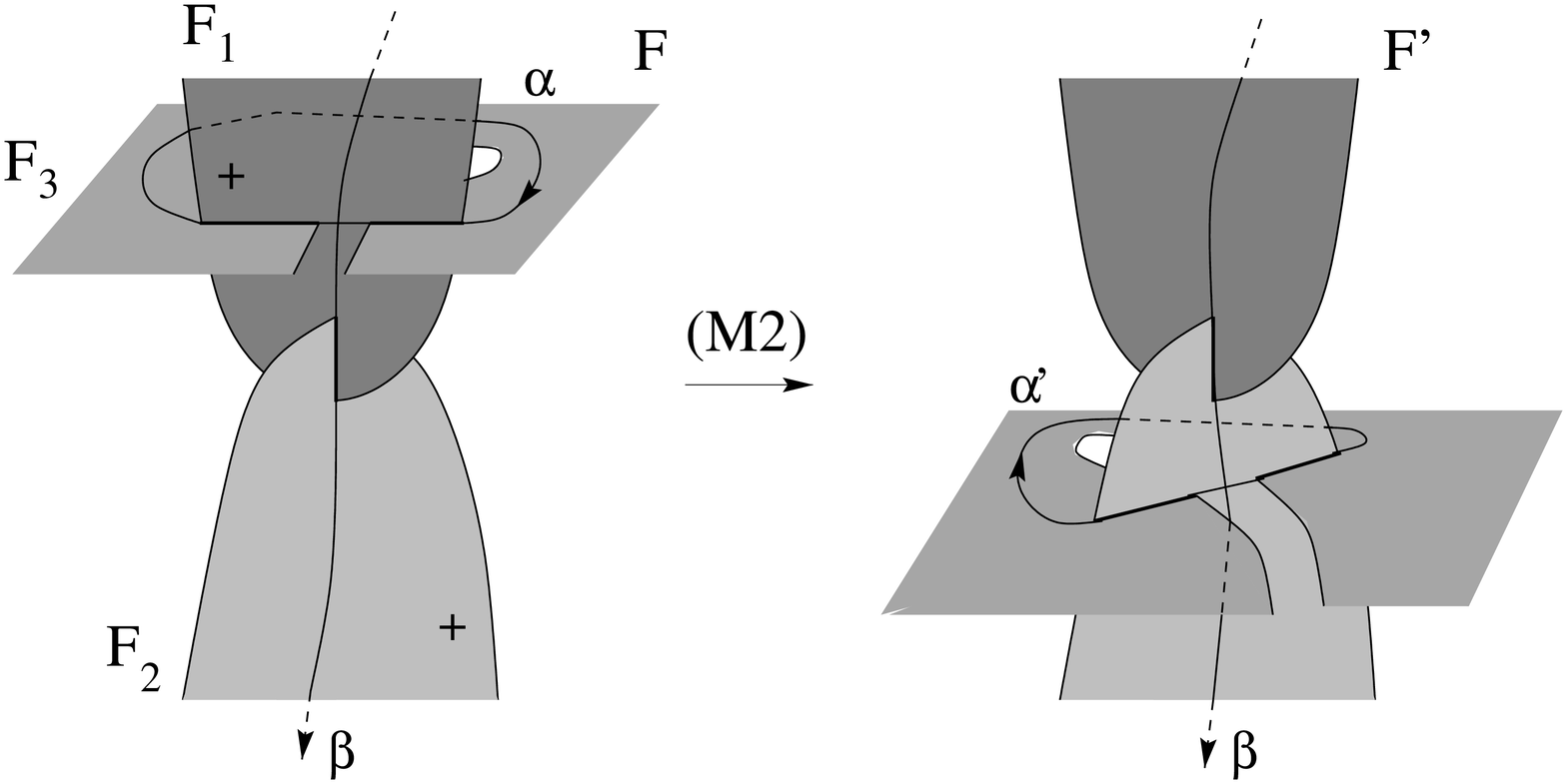,height=4.7cm}
\end{center}
\end{figure}
Here, $\Ho_1F$ and $\Ho_1F'$ are canonically isomorphic, and $\sgn(F')=\sgn(F)$. Furthermore,
$$
\chi(F'\setminus F'_i)=\cases{\chi(F\setminus F_1)-2& if $i=1$;\cr \chi(F\setminus F_2)+2& if $i=2$;\cr \chi(F\setminus F_i)& if $i>2$.}
$$
On the other hand, given any loop $\beta$ as above, we have the equalities
$$
\ell k(\alpha^\epsilon,\beta)=\cases{1& if $\epsilon(1)=+1$;\cr 0& else,}\quad\hbox{and}\quad
\ell k((\alpha')^\epsilon,\beta)=\cases{1& if $\epsilon(2)=+1$;\cr 0& else.}
$$
This implies that $(t_1-t_1^{-1})^2\det(-A_{F'})=(t_2-t_2^{-1})^2\det(-A_{F'})$. Therefore,
\begin{eqnarray*}
\Omega_{F'}(t_1,\dots,t_n)&=&\sgn(F')\prod_i(t_i-t_i^{-1})^{\chi(F'\setminus F'_i)-1}\det(-A_{F'})\\
		&=&\sgn(F)\prod_i(t_i-t_i^{-1})^{\chi(F\setminus F_i)-1}\frac{(t_2-t_2^{-1})^2}{(t_1-t_1^{-1})^2}\det(-A_{F'})\\
		&=&\Omega_F(t_1,\dots,t_n).
\end{eqnarray*}

\noindent{\bf (M3)} Consider $F'$ a C-complex obtained from $F$ by attaching a handle on $F_1$.
This time, $\Ho_1F'\simeq\Z\delta\oplus\Z\sigma\oplus\Ho_1F$, with $\delta$ and $\sigma$ as follows. 
\medskip

\begin{figure}[h]
\begin{center}
\epsfig{figure=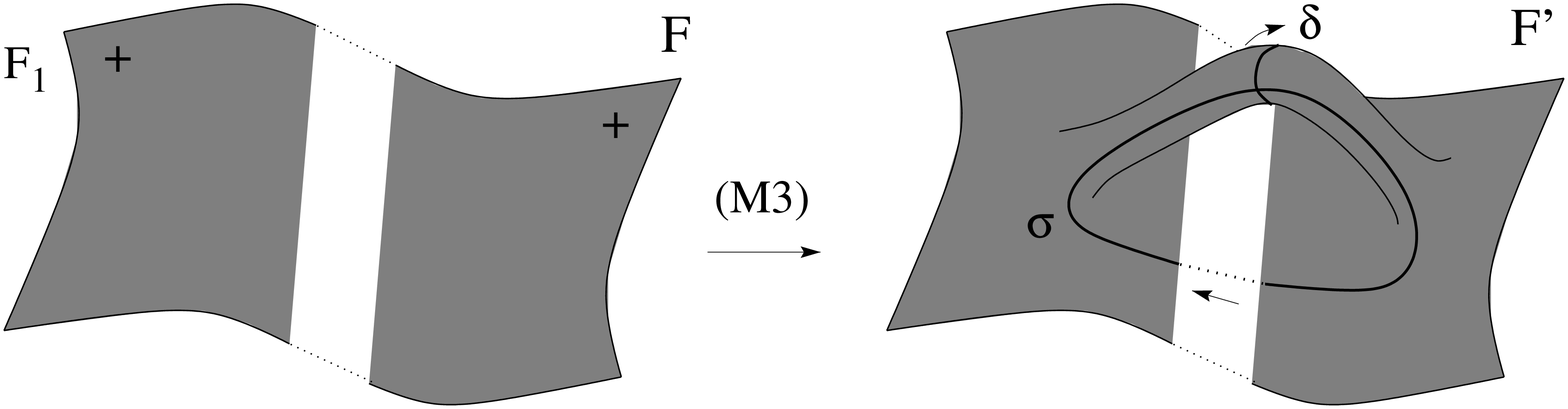,height=3cm}
\end{center}
\end{figure}
Trivially, $\sgn(F')=\sgn(F)$, and
$$
\chi(F'\setminus F'_i)=\cases{\chi(F\setminus F_i)& if $i=1$;\cr \chi(F\setminus F_i)-2& if $i>1$.}
$$
Since
$$
A_{F'}^\epsilon=\bordermatrix{&\delta&\sigma& \cr
			\delta&0&\pi(\epsilon)&0\cr
			\sigma&\pi(-\epsilon)&\ast&\ast\cr
			&0&\ast&A_F^\epsilon}\;,
\quad\hbox{with}\quad\pi(\epsilon)=\cases{-1& if $\epsilon(1)=+1$;\cr \phantom{-}0& else,}
$$
it follows
$$
-A_{F'}=\pmatrix{0&t_1\prod_{i>1}(t_i-t_i^{-1})&0\cr
		-t_1^{-1}\prod_{i>1}(t_i-t_i^{-1})&\ast&\ast\cr
		0&\ast&-A_F}.
$$
Hence, we have the equality
\begin{eqnarray*}
\Omega_{F'}(t_1,\dots,t_n)&=&\sgn(F')\prod_i(t_i-t_i^{-1})^{\chi(F'\setminus F'_i)-1}\det(-A_{F'})\\
		&=&\sgn(F)\frac{\prod_i(t_i-t_i^{-1})^{\chi(F\setminus F_i)-1}}{\prod_{i>1}(t_i-t_i^{-1})^2}\prod_{i>1}(t_i-t_i^{-1})^2\det(-A_{F})\\
		&=&\Omega_F(t_1,\dots,t_n).
\end{eqnarray*} 
This concludes the proof. \cqfd

Therefore, $\Omega_F$ is an isotopy invariant of the colored link $(L,\sigma)$; let us denote it by $\Omega_L^\sigma$.
\medskip

\noindent{\bf Lemma~5.} $\Omega_L^\sigma=\nabla_L^\sigma$.\medskip

\noindent{\it Proof.} Murakami's characterization theorem \cite{Mur} states that $\nabla^\sigma_L$ is determined uniquely by the following six relations,
where the letters $i,j$ and $k$ denote the colors of the components.
$$
\nabla^\sigma_O=\frac{1}{t_i-t_i^{-1}},\leqno(RI)
$$
where $(O,\sigma)$ denotes the trivial knot with color $i$.
$$
\nabla^\sigma_{L_+}-\nabla^\sigma_{L_-}=(t_i-t_i^{-1})\cdot \nabla^{\sigma_0}_{L_0},\leqno(RI\!I)
$$
where $L_+$, $L_-$ and $L_0$ differ by the following local operation.
\medskip

\begin{figure}[h]
\begin{center}
\epsfig{figure=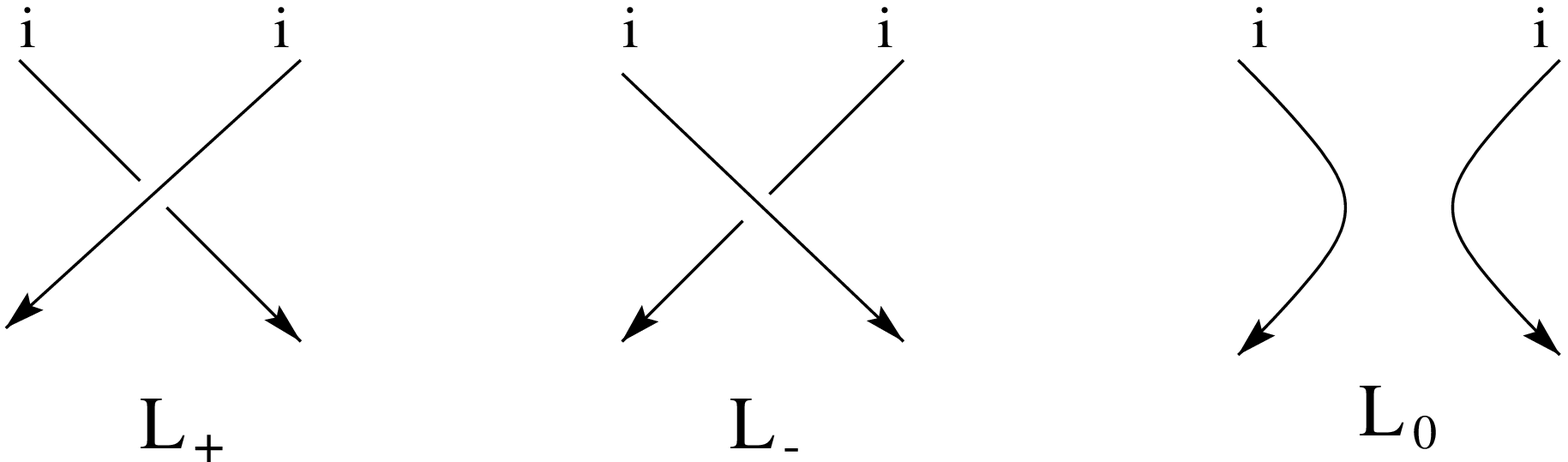,height=1.8cm}
\end{center}
\end{figure}
\vskip-0.8cm
$$
\nabla^\sigma_{L\sqcup O}=0,\leqno(RI\!I\!I)
$$
where $L\sqcup O$ denotes the disjoint union of $L$ and a trivial knot.
$$
\nabla^\sigma_{L_{++}}+\nabla^\sigma_{L_{--}}=(t_it_j+t_i^{-1}t_j^{-1})\cdot \nabla^\sigma_{L_{00}},\leqno(RIV)
$$
where $L_{++}$, $L_{--}$ and $L_{00}$ differ by the following local operation.
\medskip

\begin{figure}[h]
\begin{center}
\epsfig{figure=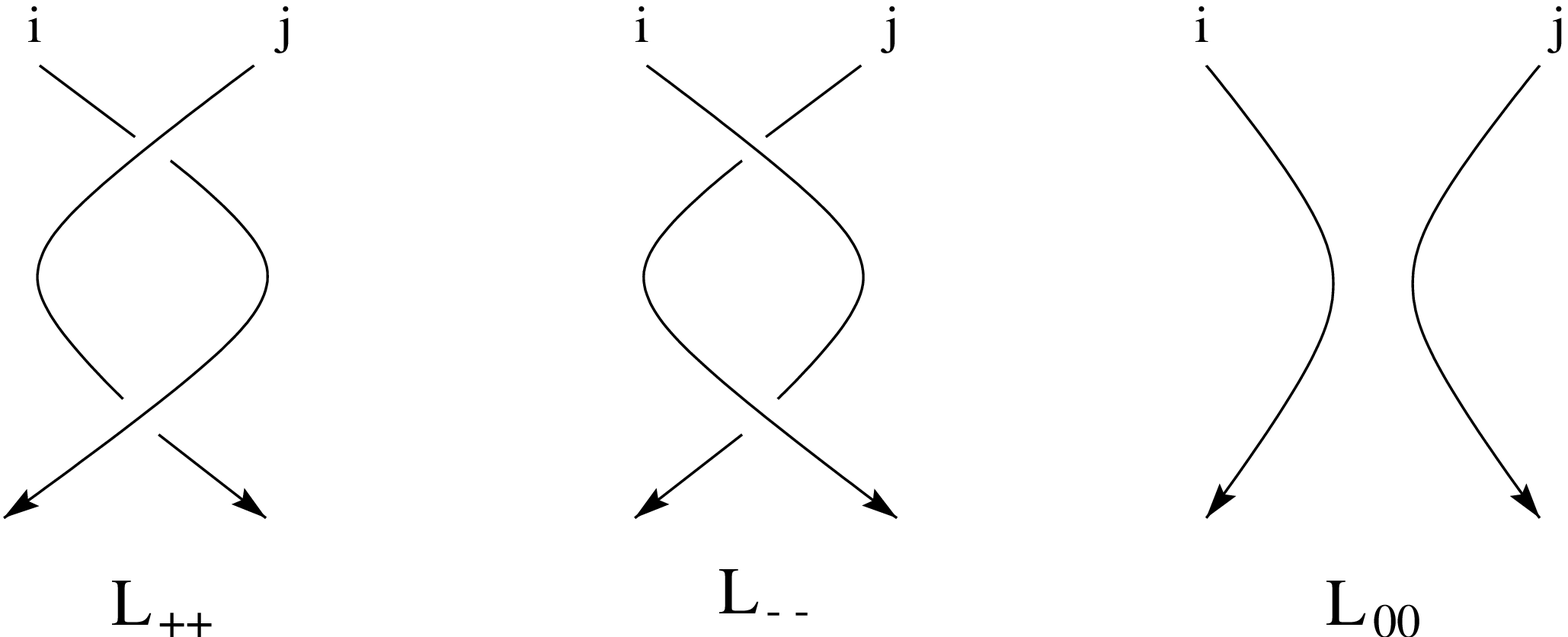,height=2.7cm}
\end{center}
\end{figure}
\vskip-0.8cm
$$
\nabla^{\sigma'}_{L'}=(t_i-t_i^{-1})\cdot\nabla^\sigma_L,\leqno(RV)
$$
where $L'$ is obtained from $L$ by the local operation given above.
\begin{figure}[h]
\begin{center}
\epsfig{figure=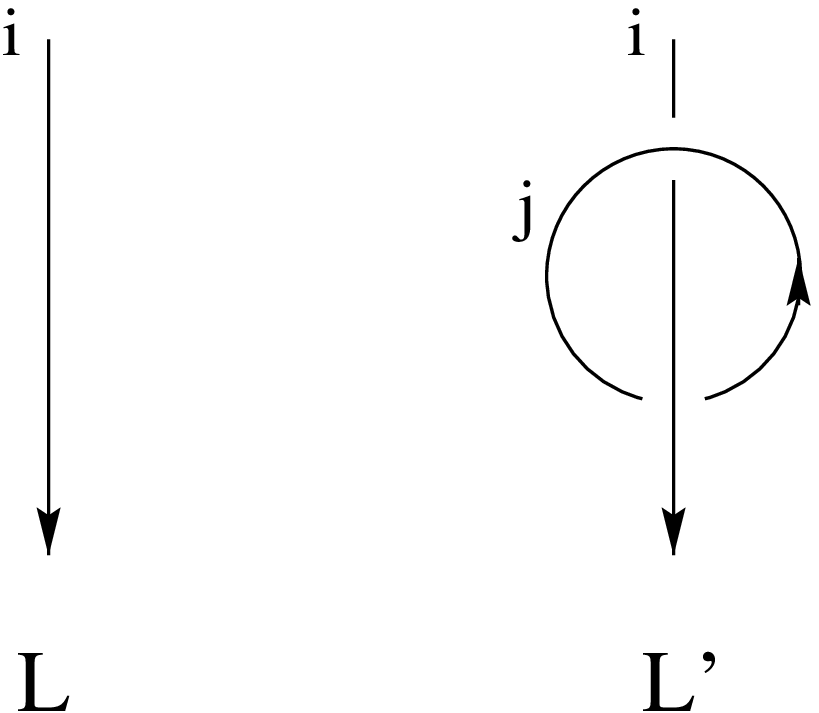,height=2.4cm}
\end{center}
\end{figure}
$$
\matrix{(t_i+t_i^{-1})(t_j-t_j^{-1})\nabla^\sigma_{L(1)}-(t_j-t_j^{-1})(t_k+t_k^{-1})\nabla^\sigma_{L(2)}\qquad\qquad\qquad\qquad\cr
-(t_i^{-1}t_k-t_it_k^{-1})(\nabla^\sigma_{L(3)}+\nabla^\sigma_{L(4)})+(t_i^{-1}t_jt_k-t_it_j^{-1}t_k^{-1})(t_k+t_k^{-1})\nabla^\sigma_{L(5)}\;\cr
-(t_i+t_i^{-1})(t_it_jt_k^{-1}-t_i^{-1}t_j^{-1}t_k)\nabla^\sigma_{L(6)}-(t_i^{-2}t_k^2-t_i^2t_k^{-2})\nabla^\sigma_{L(7)}=0\;,\qquad\;}\leqno(RV\!I)
$$
where $L(1)$, $L(2)$, $L(3)$, $L(4)$, $L(5)$, $L(6)$ and $L(7)$ differ by the following local operation.
\medskip

\begin{figure}[h]
\begin{center}
\epsfig{figure=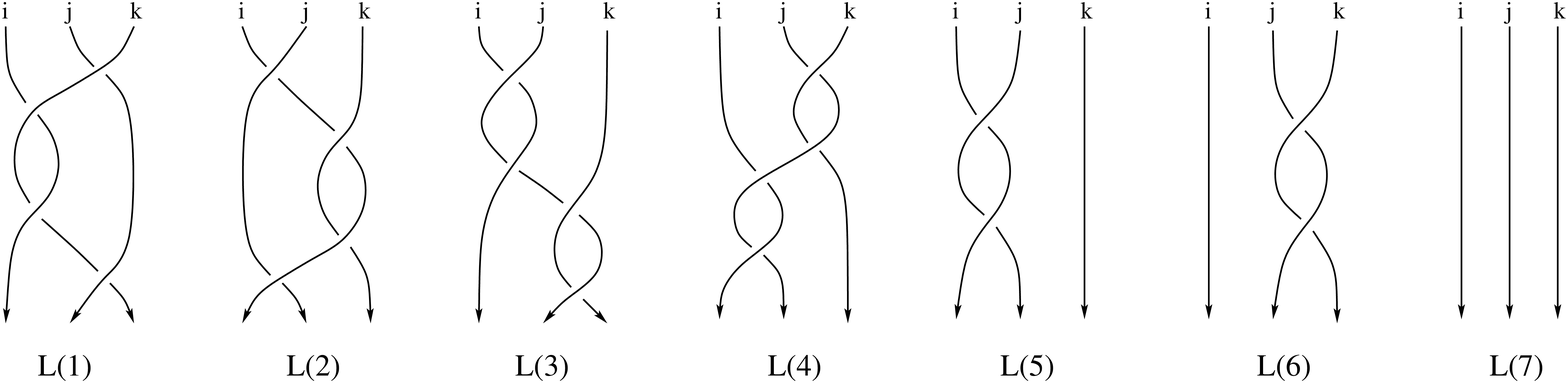,height=3cm}
\end{center}
\end{figure}
By this theorem, we just need to show that $\Omega_L^\sigma$ satisfies the relations $RI$ to $RV\!I$.\medskip

\noindent{\bf (RI)} Since a $2$-disc is a C-complex for the trivial knot, we have $\Omega^\sigma_O(t_i)=\frac{1}{t_i-t_i^{-1}}$.\medskip

\noindent{\bf (RII)} Let $F$ be a C-complex for $L_0$; C-complexes $F_+$ and $F_-$ for $L_+$ and $L_-$ are obtained from
$F$ as follows.\medskip

\begin{figure}[h]
\begin{center}
\epsfig{figure=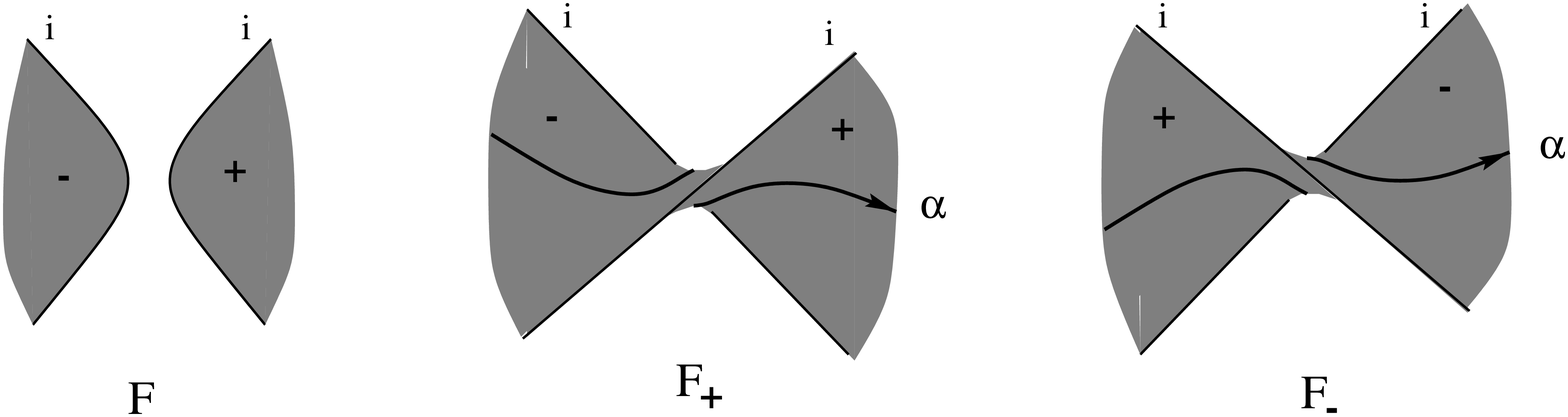,height=2.8cm}
\end{center}
\end{figure}
Clearly, $\Ho_1F_+\simeq\Ho_1F_-\simeq\Z\alpha\oplus\Ho_1F$, and for all $\epsilon$,
$$
A_{F_-}^\epsilon=\pmatrix{\ast&\ast\cr \ast& A_{F}^\epsilon}=A_{F_+}^\epsilon+\pmatrix{1&0\cr 0&0}.
$$
Therefore,
$$
A_{F_-}=\pmatrix{\ast&\ast\cr \ast& A_{F}}=A_{F_+}+\pmatrix{\prod_j(t_j-t_j^{-1})&0\cr 0&0},
$$
giving $\det\left(-A_{F_-}\right)=\det\left(-A_{F_+}\right)-\prod_j(t_j-t_j^{-1})\det\left(-A_{F}\right)$.
Since $\sgn(F_+)=\sgn(F_-)=\sgn(F)$ and
$$
\chi(F_+-(F_+)_j)=\chi(F_--(F_-)_j)=\cases{\chi(F\setminus F_j)-1& if $j\neq i$;\cr \chi(F\setminus F_i)& if $j=i$,}
$$
we get the equality
\begin{eqnarray*}
\Omega_{L_+}^\sigma-\Omega_{L_-}^\sigma&=&\sgn(F)\prod_j(t_j-t_j^{-1})^{\chi(F\setminus F_j)-2}(t_i-t_i^{-1})\left(
\det\left(-A_{F_+}\right)-\det\left(-A_{F_-}\right)\right)\\
	&=&(t_i-t_i^{-1})\cdot\sgn(F)\prod_j(t_j-t_j^{-1})^{\chi(F\setminus F_j)-1}\det\left(-A_{F}\right)\\
	&=&(t_i-t_i^{-1})\cdot\Omega_{L_0}^{\sigma_0}.
\end{eqnarray*}

\noindent{\bf (RIII)} Choose $F$ a C-complex for $L$. A C-complex $F'$ for $L\sqcup O$ is obtained from $F\sqcup D^2$ either
by the move $M1$ (if $O$ is the only component of its color) or by connecting $F_i$ and $D^2$ with a handle (if $O$ is of color $i$).
In both cases, some $1$-cycle
$x$ is created, such that $\ell k(x^\epsilon,y)=0$ for all $\epsilon$ and all $y\in\Ho_1F'$. Therefore,
$\Omega^\sigma_{L\sqcup O}$ vanishes.\medskip

\noindent{\bf (RIV)} If the colors $i$ and $j$ are equal, this relation follows easily from $RII$. Thus, we can assume
that $i=1$ and $j=2$. Given $F$ a C-complex for $L_{00}$, C-complexes $F_{++}$ for $L_{++}$ and $F_{--}$ for $L_{--}$ are obtained
as follows. 

\begin{figure}[h]
\begin{center}
\epsfig{figure=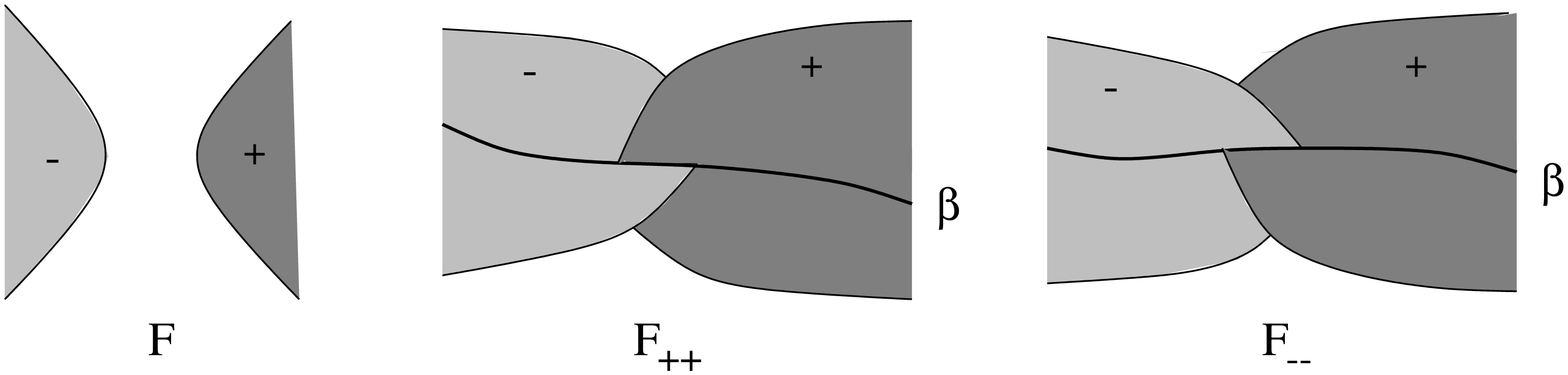,height=2.5cm}
\end{center}
\end{figure}
This time, $\Ho_1F_{++}\simeq\Ho_1F_{--}\simeq\Z\beta\oplus\Ho_1F$, and for all $\epsilon$,
$A_{F_{--}}^\epsilon=\pmatrix{\ast&\ast\cr \ast& A_{F}^\epsilon}$; on the other hand
$$
A_{F_{--}}^\epsilon=\cases{A_{F_{++}}^\epsilon& if $(\epsilon(1),\epsilon(2))=(-1,+1)$ or $(+1,-1)$;\cr
			A_{F_{++}}^\epsilon+\pmatrix{1&0\cr 0&0}& if $(\epsilon(1),\epsilon(2))=(-1,-1)$ or $(+1,+1)$.}
$$
Therefore, $\det\left(-A_{F_{--}}\right)=\det\left(-A_{F_+}\right)-
(t_1t_2+t_1^{-1}t_2^{-1})\prod_{i>2}(t_i-t_i^{-1})\det\left(-A_{F}\right)$.
Since $\sgn(F_{++})=\sgn(F)=-\sgn(F_{--})$ and
$$
\chi(F_{++}-(F_{++})_i)=\chi(F_{--}-(F_{--})_i)=\cases{\chi(F\setminus F_i)& if $i=1,2$;\cr \chi(F\setminus F_i)-1& if $i>2$,}
$$
we get the equality
\begin{eqnarray*}
\Omega_{L_{++}}^\sigma+\Omega_{L_{--}}^\sigma&=&\sgn(F)\prod_i(t_i-t_i^{-1})^{\chi(F_{++}-(F_{++})_i)-1}\left(
\det\left(-A_{F_{++}}\right)-\det\left(-A_{F_{--}}\right)\right)\\
	&=&(t_1t_2+t_1^{-1}t_2^{-1})\cdot\sgn(F)\prod_i(t_i-t_i^{-1})^{\chi(F\setminus F_i)-1}\det\left(-A_{F}\right)\\
	&=&(t_1t_2+t_1^{-1}t_2^{-1})\cdot\Omega_{L_{00}}^\sigma.
\end{eqnarray*}

\noindent{\bf (RV)} Once again, if $i=j$, this relation follows from $RII$. So, let us suppose that $i=1$ and $j=2$.
A C-complex $F'$ for $L'$ is obtained from a C-complex $F$ for $L$ by adding a disc with a positive clasp.
From the equalities $\Ho_1F'=\Ho_1F$, $\sgn(F')=\sgn(F)$ and 
$$
\chi(F'\setminus F'_i)=\cases{\chi(F\setminus F_i)+1& if $i=1$;\cr \chi(F\setminus F_i)& if $i>1$,}
$$
it follows
\begin{eqnarray*}
\Omega_{L'}^{\sigma'}&=&\sgn(F')\prod_i(t_i-t_i^{-1})^{\chi(F'\setminus F'_i)-1}\det\left(-A_{F'}\right)\\
	&=&(t_1-t_1^{-1})\cdot\sgn(F)\prod_i(t_i-t_i^{-1})^{\chi(F\setminus F_i)-1}\det\left(-A_F\right)\\
	&=&(t_1-t_1^{-1})\cdot\Omega_L^\sigma.
\end{eqnarray*}

\noindent{\bf (RVI)} Let us suppose that $i=1$, $j=2$ and $k=3$. C-complexes $F(i)$ for $L(i)$ $(i=1,\dots,6)$ are obtained from a given
C-complex $F$ for $L(7)$ as follows.
\medskip

\begin{figure}[h]
\begin{center}
\epsfig{figure=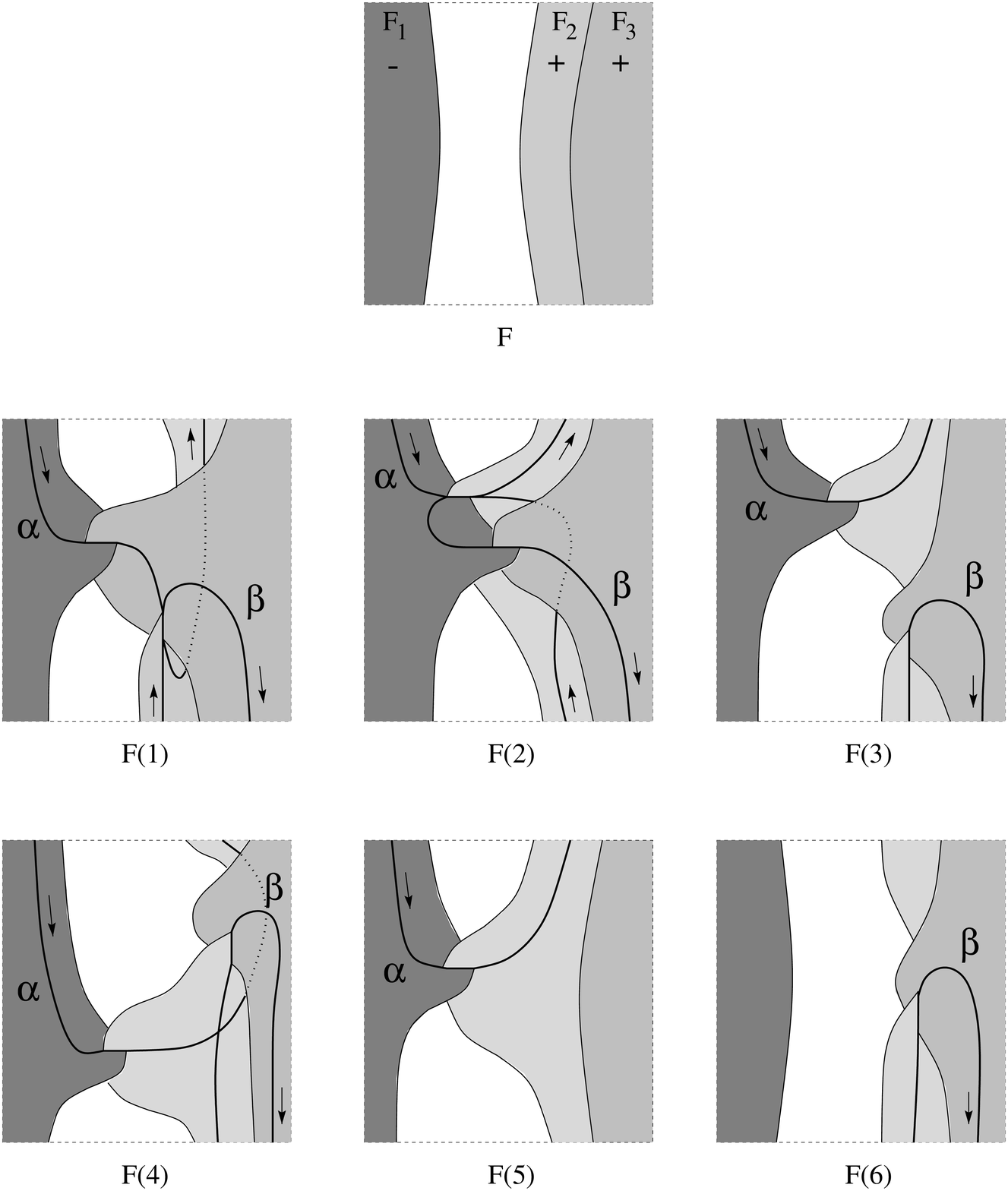,height=10cm}
\end{center}
\end{figure}
Clearly, $\Ho_1F(1)\simeq\Ho_1F(2)\simeq\Ho_1F(3)\simeq\Ho_1F(4)\simeq\Z\alpha\oplus\Z\beta\oplus\Ho_1F$, and
$\Ho_1F(5)\simeq\Z\alpha\oplus\Ho_1F$, $\Ho_1F(6)\simeq\Z\beta\oplus\Ho_1F$. By the equalities
$$
\sgn(F(i))=\sgn(F)\quad\hbox{for all}\;i\,,$$
and
\begin{eqnarray*}
\chi(F(1)\setminus F(1)_i)&=&\cases{\chi(F\setminus F_i)-1& if $i=1,2$;\cr \chi(F\setminus F_3)& if $i=3$;\cr \chi(F\setminus F_i)-2& if $i>3$,}\\
\chi(F(2)\setminus F(2)_i)&=&\cases{\chi(F\setminus F_1)& if $i=1$;\cr \chi(F\setminus F_i)-1& if $i=2,3$;\cr \chi(F\setminus F_i)-2& if $i>3$,}\\
\chi(F(3)\setminus F(3)_i)&=&\cases{\chi(F\setminus F_i)-1& if $i=1,3$;\cr \chi(F\setminus F_2)& if $i=2$;\cr \chi(F\setminus F_i)-2& if $i>3$,}=\chi(F(4)\setminus F(4)_i),\\
\chi(F(5)\setminus F(5)_i)&=&\cases{\chi(F\setminus F_i)& if $i=1,2$;\cr \chi(F\setminus F_i)-1& if $i>2$,}\\
\chi(F(6)\setminus F(6)_i)&=&\cases{\chi(F\setminus F_i)& if $i=2,3$;\cr \chi(F\setminus F_i)-1& if $i=1$ or $i>2$,}
\end{eqnarray*}
we need to check the following equation, where $T_i$ denotes $(t_i-t_i^{-1})$ and $T=\prod_{i>3}(t_i-t_i^{-1})$:
\begin{eqnarray*}
& &(t_1+t_1^{-1})\left(\frac{\det(-A_{F(1)})}{T_1\,T^2}-
(t_1t_2t_3^{-1}-t_1^{-1}t_2^{-1}t_3)\frac{\det(-A_{F(6)})}{T_1\, T}\right)\\
&-&(t_3+t_3^{-1})\left(\frac{\det(-A_{F(2)})}{T_3\, T^2}-
(t_1^{-1}t_2t_3-t_1t_2^{-1}t_3^{-1})\frac{\det(-A_{F(5)})}{T_3\, T}\right)\\
&=&(t_1^{-1}t_3-t_1t_3^{-1})\left(\frac{\det(-A_{F(3)})+\det(-A_{F(4)})}{T_1\,T_3\, T^2}+
(t_1^{-1}t_3+t_1t_3^{-1})\det(-A_F)\right).
\end{eqnarray*}
Without loss of generality, it may be assumed that $F$ has the following form near the place of the skein relation.

\begin{figure}[h]
\begin{center}
\epsfig{figure=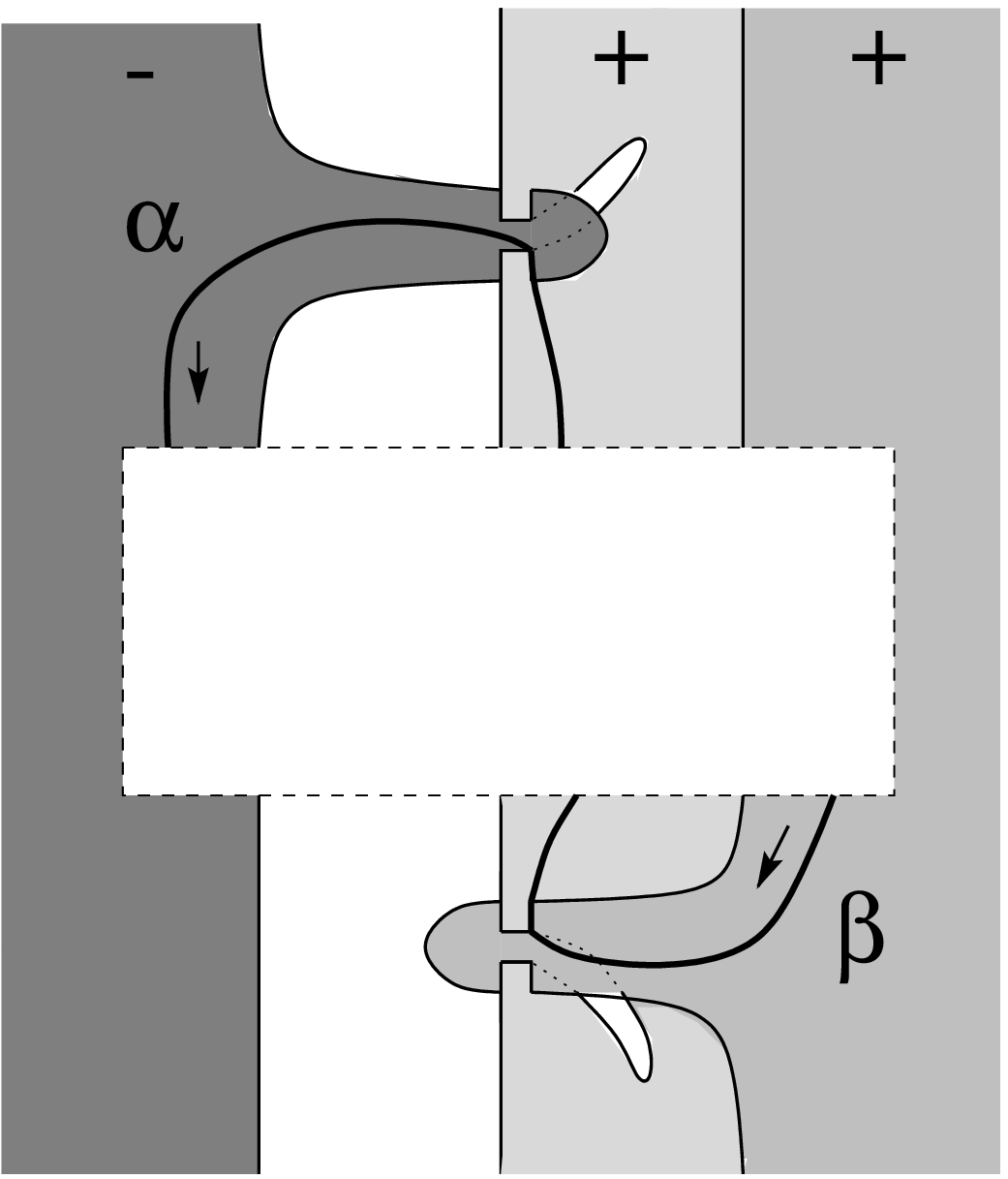,height=4.5cm}
\end{center}
\end{figure}
It is then easy to compute the matrices $A_{F(i)}$:
\begin{eqnarray*}
-A_{F(1)}&=& \bordermatrix{&\alpha&\beta& \cr 
	\alpha&(t_1t_2t_3^{-1}-t_1^{-1}t_2^{-1}t_3)\, T& -T_1\,t_2\,t_3\,T& v_\alpha\,T_3\,T\cr
	\beta&-T_1\,t_2^{-1}t_3^{-1}\,T&0&v_\beta\,T_1\,T\cr
	&w_\alpha\,T_3\,T&w_\beta\,T_1\,T&A_F},\\
-A_{F(2)}&=& \bordermatrix{&\alpha&\beta& \cr 
	\alpha&0& -t_1\,t_2\,T_3\,T& v_\alpha\,T_3\,T\cr
	\beta&-t_1^{-1}t_2^{-1}\,T_3\,T&(t_1^{-1}t_2t_3-t_1t_2^{-1}t_3^{-1})&v_\beta\,T_1\,T\cr
	&w_\alpha\,T_3\,T&w_\beta\,T_1\,T&A_F},\\
-A_{F(3)}&=& \bordermatrix{&\alpha&\beta& \cr 
	\alpha&0&0& v_\alpha\,T_3\,T\cr
	\beta&0&0&v_\beta\,T_1\,T\cr
	&w_\alpha\,T_3\,T&w_\beta\,T_1\,T&A_F},\\
-A_{F(4)}&=& \bordermatrix{&\alpha&\beta& \cr 
	\alpha&0&T_1\,t_2\,T_3\,T&v_\alpha\,T_3\,T\cr
	\beta&-T_1\,t_2^{-1}\,T_3\,T&0&v_\beta\,T_1\,T\cr
	&w_\alpha\,T_3\,T&w_\beta\,T_1\,T&A_F},\\
-A_{F(5)}&=& \bordermatrix{&\alpha& \cr
			\alpha&0& v_\alpha\,T_3\,T\cr
			&w_\alpha\,T_3\,T&A_F},\\
-A_{F(6)}&=&  \bordermatrix{&\beta& \cr
			\beta&0& v_\beta\,T_1\,T\cr
			&w_\beta\,T_1\, T&A_F}.\\
\end{eqnarray*}
The equality can then be checked by direct computation. This concludes Lemma~5, and the theorem is proved. \cqfd

\section{Applications}

In his fundamental paper \cite{Con}, Conway states various properties of his potential function.
We show now that these properties follow immediately from our construction.
\medskip

\noindent{\bf Proposition~1.} $\nabla^\sigma_L(t_1^{-1},\dots,t_n^{-1})=(-1)^\mu\,\nabla_L^{\sigma}(t_1,\dots,t_n).$
\medskip

\noindent{\it Proof.} By the theorem, we easily get $\nabla_L(t_1^{-1},\dots,t_\mu^{-1})=(-1)^\nu\nabla_L(t_1,\dots,t_\mu)$, with
$$
\nu=\sum_{i=1}^n(\chi(F\setminus F_i)-1)+n\cdot\rk\Ho_1F.
$$
For $i=1,\dots,n$, let us denote by $\mu_i$ the number of components of $L$ with color $i$.
Computing modulo $2$, we have
\begin{eqnarray*}
\nu&\equiv&\sum_{i=1}^n\chi(F\setminus F_i)+n\cdot\left(\rk\Ho_1F+1\right)\\
	&\equiv&\sum_{i=1}^n(\chi(F)+\overbrace{\chi(F_i)}^{\equiv \mu_i}+\chi(F_i\cap(F\setminus F_i)))+n\cdot\chi(F)\\
	&\equiv&\sum_{i=1}^n \mu_i + \sum_{i=1}^n\chi(F_i\cap(F\setminus F_i))=\mu+2\cdot\#\{\hbox{clasps}\}\equiv\mu,
\end{eqnarray*}
and the proposition is proved. \cqfd

\noindent {\bf Proposition~2.} $\nabla_L^\sigma(-t_1,\dots,-t_n)=(-1)^\mu\,\nabla_L^{\sigma}(t_1,\dots,t_n)$.
\medskip

\noindent{\it Proof.} Exactly as Proposition~1. \cqfd

\noindent {\bf Proposition~3.} {\it Let $(L,\sigma)$ be a colored link with $\#\sigma^{-1}(1)=\mu_1$, and let
$(L',\sigma)$ be the same link with the opposite orientation on the sublink of color $1$. Then, we have the equality}
$$
\nabla_{L'}^\sigma(t_1,\dots,t_n)=(-1)^{\mu_1}\,\nabla_L^{\sigma}(t_1^{-1},t_2,\dots,t_n).
$$
\noindent{\it Proof.} If $F=F_1\cup\dots\cup F_n$ is a C-complex for $L$, a C-complex for $L'$ is given by $F'=(-F_1)\cup F_2\cup\dots\cup F_n$.
Since $\sgn(F')=(-1)^{\#\{\hbox{\scriptsize{clasps of }}F_1\}}\cdot\sgn(F)=(-1)^{\chi(F_1\cap(F\setminus F_1))}\cdot\sgn(F)$ and
$A^\epsilon_{F'}=A^{\epsilon'}_F$ with $\epsilon'(1)=-\epsilon(1)$ and $\epsilon'(i)=\epsilon(i)$ for $i>1$, we have
$$
\nabla_{L'}(t_1,\dots,t_\mu)=(-1)^\nu\,\nabla_L(t_1^{-1},t_2,\dots,t_\mu),
$$
where
\begin{eqnarray*}
\nu&=&\chi(F_1\cap(F\setminus F_1))+\chi(F\setminus F_1)-1+\rk\Ho_1F\\
&\equiv&\chi(F)-\chi(F\setminus F_1)+\chi(F_1\cap(F\setminus F_1))=\chi(F_1)\equiv\mu_1.
\end{eqnarray*}
This settles the proof. \cqfd

\noindent{\bf Corollary~1.} $\nabla^\sigma_{-L}(t_1,\dots,t_n)=(-1)^\mu\,\nabla_L^{\sigma}(t_1^{-1},\dots,t_n^{-1}).$ \cqfd

\noindent{\bf Corollary~2.} $\nabla^\sigma_{-L}=\nabla_L^{\sigma}.$ \cqfd

\noindent{\bf Proposition~4.} {\it $\nabla_{m(L)}^\sigma=(-1)^{\mu+1}\,\nabla_L^\sigma$, where $m(L)$ denotes the mirror image of $L$.}
\medskip

\noindent{\it Proof.} The mirror image $m(F)$ of a C-complex
$F$ for $L$ provides a C-complex for $m(L)$. The sign of every clasp and of every linking number is changed, giving
$A^\epsilon_{m(F)}=-A^\epsilon_F$ for all $\epsilon$. Hence
$$
\nabla_{m(L)}=(-1)^{\#\{\hbox{\scriptsize{clasps of }}F\}}\cdot(-1)^{\hbox{\scriptsize{rk}}\,\hbox{\scriptsize{H}}_1F}\cdot\nabla_L.
$$
The equality $\chi(F)=\sum_{i=1}^n\chi(F_i)-\#\{\hbox{clasps of $F$}\}\,$ yields
$$
\#\{\hbox{clasps of $F$}\}+\rk\Ho_1F=\rk\Ho_1F-\chi(F)+\sum_{i=1}^\mu\overbrace{\chi(F_i)}^{\equiv \mu_i}\equiv \mu+1,
$$
and the proposition is proved. \cqfd

Let $(L',\sigma')$ and $(L'',\sigma'')$ be two colored links; let us suppose that two components $L_i'$ of $L'$ and $L_j''$
of $L''$ have the same color $\sigma'(i)=\sigma''(j)=\alpha$. Then, there is a well-defined colored link given by the connected sum of $L'$ and
$L''$ along $L_i'$ and $L_j''$; we will denote it by $(L'\#L'',\sigma'\#\sigma'')$. The following multiplicativity formula was also stated
by Conway. As far as I know, there is no published proof.
\medskip

\noindent{\bf Proposition~5.} $\nabla^{\sigma'\#\sigma''}_{L'\#L''}=\nabla_{L'}^{\sigma'}\cdot\nabla_{L''}^{\sigma''}
\cdot(t_\alpha-t_\alpha^{-1}).$
\medskip

\noindent{\it Proof.} Given C-complexes $F'$ for $(L',\sigma')$ and $F''$ for $(L'',\sigma'')$, a C-complex for the connected sum 
$(L'\#L'',\sigma'\#\sigma'')$ is given by $F=F'\#F''$. Clearly, $A^\epsilon_F=A^\epsilon_{F'}\oplus A^\epsilon_{F''}$ for all $\epsilon$,
giving
$$
A_F=\left(\prod_{j\neq\alpha}\left(t_j''-(t_j'')^{-1}\right)\cdot A_{F'}\right)\oplus
\left(\prod_{i\neq\alpha}\left(t_i'-(t_i')^{-1}\right)\cdot A_{F''}\right).
$$
Since $\sgn(F)=\sgn(F')\cdot\sgn(F'')$ and
\begin{eqnarray*}
\chi(F\setminus F_i')&=&\chi(F'\setminus F_i')+\chi(F'')-1\quad\forall\;i\neq\alpha,\\
\chi(F\setminus F_j'')&=&\chi(F')+\chi(F''\setminus F_j'')-1\quad\forall\;j\neq\alpha,\\
\chi(F\setminus F_\alpha)&=&\chi(F'\setminus F'_\alpha)+\chi(F''\setminus F''_\alpha),
\end{eqnarray*}
we get the result. \cqfd

For the sake of completeness, let us mention without proof two additional properties of Conway's potential function.
The first one is another skein relation announced by Conway, that can be proved using our construction.
\medskip

\noindent{\bf Proposition~6.} {\it $\nabla_{L_1}^\sigma+\nabla_{L_2}^\sigma=\nabla_{L_3}^\sigma+\nabla_{L_4}^\sigma$, where
$L_1,L_2,L_3$ and $L_4$ are identical except within a ball where they are related as illustrated below.}

\begin{figure}[h]
   \begin{center}
     \epsfig{figure=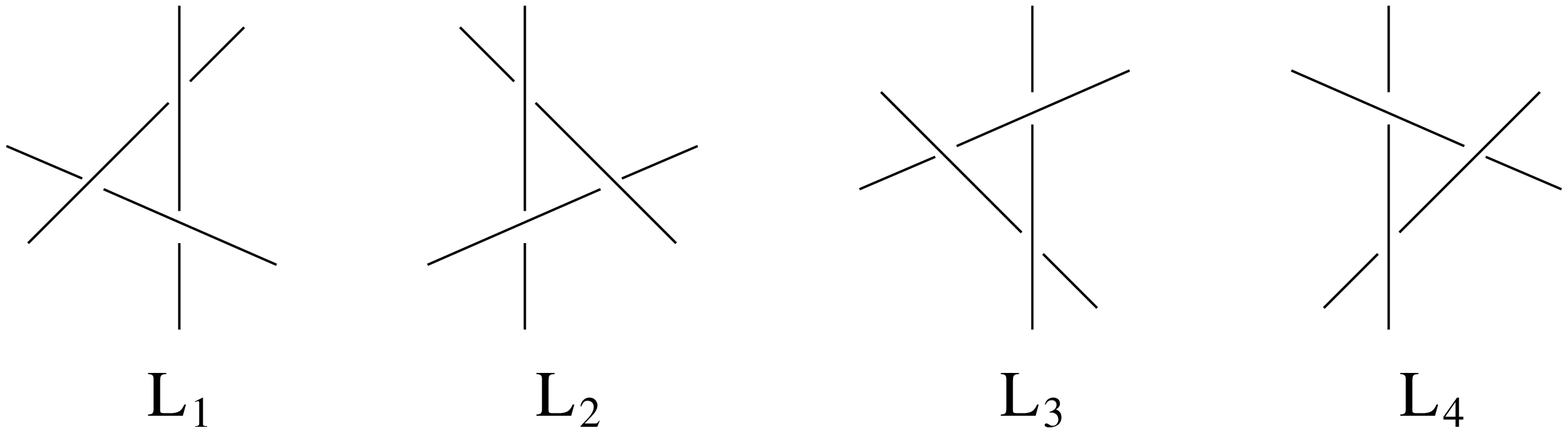,height=2.3cm}
   \end{center}
\end{figure}
\vskip-1cm\cqfd\vskip0.7cm

The second one is the translation of the Torres formula \cite[Theorem 3]{Torres} from the Alexander polynomial to the Conway potential function.
This formula can be recovered from our theorem, but the proof is a little tedious.
\medskip

\noindent{\bf Proposition~7.} {\it $\nabla_L(1,t_2,\dots,t_\mu)=(t_2^{\ell_{12}}\cdots t_\mu^{\ell_{1\mu}}-
t_2^{-\ell_{12}}\cdots t_\mu^{-\ell_{1\mu}})\cdot\nabla_{L\setminus L_1}(t_2,\dots,t_\mu)$, where $\ell_{1i}$ denotes the linking number
$\ell k(L_1,L_i)$. \cqfd }

\nonumsection{Acknowledgements}
The author wishes to express his thanks to Daryl Cooper and to Mathieu Baillif. Above all, it is a pleasure to thank Jerry Levine
for his kindness and disponibility at Brandeis University.								

\nonumsection{References}

\end{document}